\documentclass[12pt,a4paper,twoside]{article}

\usepackage[ngerman, english]{babel}
\usepackage[T1]{fontenc}
\usepackage{amsmath,amssymb, amsthm}
\usepackage[DIV14,BCOR2cm]{typearea}
\usepackage{enumerate}
\usepackage[arrow, matrix, curve]{xy}
\usepackage{mathdots}
\usepackage{verbatim}
\usepackage{mathabx}
\usepackage{blkarray}
\usepackage{multirow}
\usepackage{bbm}
\usepackage{hyperref}

\newtheoremstyle{thm}%
        {10pt}
        {10pt}
        {\itshape}
        {}
        {\bfseries}
        {}
        {2mm}
        {}
\newtheoremstyle{def}%
        {10pt}
        {10pt}
        {\normalfont}
        {}
        {\bfseries}
        {}
        {2mm}
        {}
\theoremstyle{def}
\newtheorem{Definition}{Definition}[section]

\newtheorem*{Remark}{Remark}
\newtheorem*{Proof}{Proof}
\theoremstyle{thm}
\newtheorem{Theorem}[Definition]{Theorem}
\newtheorem{Lemma}[Definition]{Lemma}

\newtheorem{Corollary}[Definition]{Corollary}
\numberwithin{equation}{section}
\title{\bf Structure of Hyperbolic Unitary Groups II: Classification of E-normal Subgroups}
\author{{Raimund Preusser} 
\\ E-mail: preusser@math.uni-bielefeld.de}
\date{}
\begin{document}
\maketitle
\noindent
{\bf Abstract.} This paper proves the sandwich classification conjecture for subgroups of an even dimensional hyperbolic unitary group $U_{2n}(R,\Lambda)$ which are normalized by the elementary subgroup $EU_{2n}(R,\Lambda)$, under the condition that $R$ is a quasi-finite ring with involution, i.e a direct limit of module finite rings with involution, and $n\geq 3$.\\
\\
{\bf 2010 Mathematics Subject Classification:} 20G35, 20H25\\
\\
{\bf Keywords:} unitary groups, E-normal subgroups, sandwich classification

\section{Introduction}
This paper is a successor of the paper \cite{bak-vavilov} by A. Bak and N. Vavilov.
The main result is the following: If $(R,\Lambda)$ is a form ring such that $R$ is quasi-finite and $H$ is a subgroup of the hyperbolic unitary group $U_{2n}(R,\Lambda)$ where $n\geq 3$, then 
\begin{align*}  &H \text{ is normalized by the elementary subgroup } EU_{2n}(R,\Lambda)\text{ of } U_{2n}(R,\Lambda)\\ 
\Leftrightarrow~&\exists!\text{ form ideal }(I,\Gamma):EU_{2n}((R,\Lambda),(I,\Gamma))\subseteq H\subseteq CU_{2n}((R,\Lambda),(I,\Gamma))\tag{1.1}
\end{align*}
where $EU_{2n}((R,\Lambda),(I,\Gamma))$ denotes the relative elementary subgroup of level $(I,\Gamma)$ and $CU_{2n}((R,\Lambda),(I,\Gamma))$ denotes the full congruence subgroup of level $(I,\Gamma)$. This result extends the range of validity of previous results. If $R$ has finite Bass-Serre dimension $d$ (cf.\hspace{-0.15cm} \cite{bak}) then the result was proved already in \cite{bak_2} provided $n\geq sup(d+2,3)$ and if $R$ is commutative, it was proved recently in \cite{you-zhou}. An incorrect proof, which can be repaired when $2$ is invertible in $R$, was given in \cite{vaserstein-you}. \\

The paper is organized as follows. In section 2 we recall some standard notation which will be used throughout the paper. In section 3 we recall the definitions of 
the hyperbolic unitary group and some important subgroups. In section 4 we prove the main result (1.1), first for certain almost commutative rings and then for quasi-finite rings. \\

The current paper formed a part of my doctoral dissertation. I would like to thank my advisor Anthony Bak for his guidance during the preparation of my dissertation and in particular for making me aware of the theory of model unitary groups, which is used in section 4.
\section{Notation} 
Let $G$ be a group and $H,K$ be subsets of $G$. The subgroup of $G$ generated by $H$ is denoted by $\langle H\rangle$. If $g,h\in G$, let $^hg:=hgh^{-1}$, $g^h:=h^{-1}gh$ and $[g,h]:=ghg^{-1}h^{-1}$. Set $^KH:=\langle \{^kh|h\in H, k\in K\}\rangle$ and $H^K:=\langle \{h^k|h\in H, k\in K\}\rangle$. Analogously define $[H,K]$ and $HK$. Instead of $^K\{g\}$ we write $^Kg$ (analogously we write $g^K$ instead of $\{g\}^K$, $^gH$ instead of $^{\{g\}}H $, $[g,K]$ instead of $[\{g\},K]$ etc.).\\

In this paper, ring will always mean associative ring with $1$ such that $1\neq 0$. Ideal will mean two-sided ideal. If $R$ is a ring and $m,n\in \mathbb{N}$, then the set of all invertible elements in $R$ is denoted by $R^*$ and the set of all $m\times n$ matrices with entries in $R$ is denoted by $M_{m\times n}(R)$. If $a\in M_{m\times n}(R)$, let $a_{ij}\in R$ denote the element in the $(i,j)$'th position. Let $a^t\in M_{n\times m}(R)$ denote its transpose, thus $(a^t)_{ij}=a_{ji}$. Denote the $i$-th row of $a$ by $a_{i*}$ and the $j$-th column of $a$ by $a_{*j}$. We set $M_n(R):=M_{n\times n}(R)$. The identity matrix in $M_n(R)$ is denoted by $e$ or $e^{n\times n}$ and the matrix with a $1$ at position $(i,j)$ and zeros elsewhere is denoted by $e^{ij}$. If $a\in M_n(R)$ is invertible, the entry 
of $a^{-1}$ at position $(i,j)$ is denoted by $a'_{ij}$, the $i$-th row of $a^{-1}$ by $a'_{i*}$ and the $j$-th column of $a^{-1}$ by $a'_{*j}$. Further we denote by $^nR$ the set of all rows $v=(v_1,\dots,v_n)$ with entries in $R$ and by $R^n$ the set of all columns $u=(u_1,\dots,u_n)^t$ with entries in $R$.
\section{Bak's hyperbolic unitary groups}
In order to classify the subgroups of a general linear
group normalized by its elementary subgroup, the notion of an ideal in a ring is sufficient. Bak's dissertation \cite{bak_2} showed  that the notion of an ideal by itself
 was not sufficient to solve the analogous classification problem for unitary groups, but that a refinement
of the notion an ideal, called a form ideal, was necessary. This led naturally to a more general notion 
of unitary group, which was defined over a form ring instead of just a ring and generalized all previous concepts. We describe form rings $(R,\Lambda)$ and form ideals ideals $(I, \Gamma)$ first, then hyperbolic unitary 
groups $U_{2n}(R,\Lambda)$ over form rings $(R,\Lambda)$. For form ideals $(I,\Gamma)$, we recall the definitions 
 of the following subgroups of $U_{2n}(R,\Lambda)$; the preelementary groups $EU_{2n}(I, \Gamma)$, 
the relative elementary groups $EU_{2n}((R,\Lambda),(I,\Gamma))$, the principal congruence subgroups $U_{2n}((R,\Lambda),(I,\Gamma))$, 
and the full congruence subgroups $CU_{2n}((R,\Lambda),(I,\Gamma))$.
\Definition{Let $R$ be a ring and \begin{align*}\bar{}:R&\rightarrow R\\ r&\mapsto \overline{r}\end{align*}
an involution on $R$, i.e. $\overline{r+s}=\overline{r}+\overline{s}$, $\overline{rs}=\bar{s}\bar{r}$ and $\overline{\overline{r}}=r$ for any $r,s\in R$. Let $\lambda\in center(R)$ such that $\lambda\overline{\lambda}=1$ and set $\Lambda_{min}:=\{r-\lambda\overline{r}|r\in R\}$ and $\Lambda_{max}:=\{r\in R|r=-\lambda\overline{r}\}$. An additive subgroup $\Lambda$ of $R$ such that 
\begin{enumerate}[(1)]
\item $\Lambda_{min}\subseteq\Lambda\subseteq\Lambda_{max}$ and
\item $r\Lambda\overline{r}\subseteq\Lambda~\forall r\in R$ 
\end{enumerate}
is called a {\it form parameter}. If $\Lambda$ is a form parameter for R, the pair $(R,\Lambda)$ is called a {\it form ring}.
}
\Definition{Let $(R,\Lambda)$ be a form ring and $I$ an ideal such that $\overline{I}=I$. Set $\Gamma_{max}=I\cap\Lambda$ and $\Gamma_{min}=\{\xi-\lambda\overline{\xi}|\xi\in I\}+\langle \{\zeta\alpha\overline{\zeta}|\zeta\in I,\alpha\in\Lambda\}\rangle$. If we want to stress that $\Gamma_{max}$ (resp. $\Gamma_{min}$) belongs to $I$, we write $\Gamma_{max}^I$ (resp. $\Gamma_{min}^I$). An additive subgroup $\Gamma$ of $I$ such that 
\begin{enumerate}[(1)]
\item $\Gamma_{min}\subseteq\Gamma\subseteq\Gamma_{max}$ and  
\item $\alpha\Gamma\overline{\alpha}\subseteq\Gamma~\forall \alpha\in R$ 
\end{enumerate}
is called a {\it relative form parameter of level} $I$. If $\Gamma$ is a relative form parameter of level $I$, then $(I,\Gamma)$ is called a {\it form ideal} of $(R,\Lambda)$.
}\\

Until the end of this section let $n\in\mathbb{N}$, $(R,\Lambda)$ a form ring and $(I,\Gamma)$ a form ideal of $(R,\Lambda)$.\\
\Definition{Let $V$ be a free right $R$-module of rank $2n$ and $B=(e_1,\dots,e_n,e_{-n},$ $\dots,e_{-1})$ an ordered basis of $V$. Let $\phi_B:V \rightarrow R^{2n}$ be the module isomorphism mapping $e_i$ to the column whose $i$-th coordinate is one and all the other coordinates are zero if $1\leq i\leq n$ and the column whose $(2n+1+i)$-th coordinate is one and all the other coordinates are zero if $-n\leq i \leq -1$. In the following we will identify elements $v\in V$ with their images $\phi_B(v)\in R^{2n}$. Let $p\in M_n(R)$ be the matrix with ones on the skew diagonal and zeros elsewhere. We define the maps 
\begin{align*}
\mathbbm{f}:V\times V&\rightarrow R\\(v,w)&\mapsto \overline{v}^t\begin{pmatrix} 0 & p \\ 0 & 0 \end{pmatrix}w,
\end{align*} 
\begin{align*}
\mathbbm{h}:V\times V&\rightarrow R\\(v,w)&\mapsto \overline{v}^t\begin{pmatrix} 0 & p \\ \lambda p & 0 \end{pmatrix}w
\end{align*} 
and
\begin{align*}
\mathbbm{q}:V&\rightarrow R/\Lambda\\v&\mapsto \mathbbm{f}(v,v)+\Lambda
\end{align*} 
where $\bar v$ is obtained from $v$ by applying $~\bar{}~$ to each entry of $v$. The maps $\mathbbm{f}$, $\mathbbm{h}$ and $\mathbbm{q}$ are denoted in \cite{bak-vavilov}, page 164, by $f$, $h$ and $q$, respectively.
It is easy to check that $\mathbbm{f}(v,w)=\overline{v}_1w_{-1}+...+\overline{v}_nw_{-n}$, $\mathbbm{h}(v,w)=\overline{v}_1w_{-1}+...+\overline{v}_nw_{-n}+\lambda\overline{v}_{-n}w_n+...+\lambda\overline{v}_{-1}w_1=\mathbbm{f}(v,w)+\lambda\overline{\mathbbm{f}(w,v)}$ and $\mathbbm{q}(v)=\overline{v}_1v_{-1}+...+\overline{v}_nv_{-n}+\Lambda$ for any $v,w\in V$. For any $v\in V$, $\mathbbm{f}(v,v)$ is called the {\it length} of $v$ and is denoted by $|v|$.
}
\Definition{The subgroup $U_{2n}(R,\Lambda):=\{\sigma\in GL(V)|(\mathbbm{h}(\sigma u,\sigma v)=\mathbbm{h}(u,v)) \wedge (\mathbbm{q}(\sigma u)=\mathbbm{q}(u))~\forall u,v\in V\}$ of $GL(V)$ is called the {\it hyperbolic unitary group}. We will identify $U_{2n}(R,\Lambda)$ with its image in $GL_{2n}(R)$ under the isomorphism $GL(V)\rightarrow GL_{2n}(R)$ determined by the ordered basis $(e_1,\dots,e_n,e_{-n},\dots,e_{-1})$. 
}
\Definition{Let $\sigma\in M_n(R)$. By definition $\sigma^*$ is the matrix in $M_n(R)$ whose entry at position $(i,j)$ equals $\overline{\sigma}_{ji}$. Further we define $AH_n(R,\Lambda):=\{a\in M_n(R)|a=-\lambda a^*, a_{ii}\in \Lambda~\forall i\in\{1,...,n\}\}$.}
\Lemma{Let $(R,\Lambda)$ be a form ring, $n\in \mathbb{N}$ and $\sigma=\begin{pmatrix} a & b \\ c & d \end{pmatrix}\in GL_{2n}(R)$, where $a,b,c,d\in M_n(R)$. Then $\sigma\in U_{2n}(R,\Lambda)$ if and only if
\begin{enumerate}[(1)]
\item $\sigma^{-1}=\begin{pmatrix} pd^*p & \overline{\lambda}pb^*p \\ \lambda pc^*p & pa^*p \end{pmatrix}$ and
\item $a^*pc, b^*pd\in AH_n(R,\Lambda)$.
\end{enumerate}
}
\Proof{See \cite{bak-vavilov}, p.166.}
\Remark{\hspace{1cm}
\begin{enumerate}[(1)]
\item If $a\in M_n(R)$, then $pa^*p$ is the matrix one gets by applying the involution to each entry of $a$ and mirroring all entries on the skew diagonal.
\item In $\cite{bak_2}$, $\cite{you}$ and $\cite{you-zhou}$ the ordered basis $(e_1,...,e_n,e_{-1},...,e_{-n})$ is used and hence the matrices may look different. Let $\sigma\in GL(V)$. If the image of $\sigma$ under the isomorphism $GL(V)\rightarrow GL_{2n}(R)$ determined by the ordered basis $(e_1,...,e_n,e_{-1},$ $...,e_{-n})$ (which is used in the papers mentioned above) equals $\begin{pmatrix} a & b \\ c & d \end{pmatrix}$, where $a,b,c,d\in M_n(R)$, then the image of $\sigma$ under the isomorphism $GL(V)\rightarrow GL_{2n}(R)$ determined by the ordered basis $(e_1,...,e_n,e_{-n},...,e_{-1})$ (which is used in this paper) equals $\begin{pmatrix} a & bp \\ pc & pdp \end{pmatrix}$.
\end{enumerate}
\Definition{We define $\Omega_+:=\{1,...,n\}$, $\Omega_-:=\{-n,...,-1\}$, $\Omega:=\Omega_+\cup\Omega_-$ and 
\begin{align*}
\epsilon:\Omega&\rightarrow\{-1,1\}\\
i&\mapsto\epsilon(i):=\begin{cases}
1, &if~i\in\Omega_+,\\
-1, &if~i\in\Omega_-. 
\end{cases}
\end{align*}
\Lemma{Let $\sigma\in GL_{2n}(R)$. Then $\sigma\in U_{2n}(R,\Lambda)$ if and only if
\begin{enumerate}[(1)]
\item $\sigma'_{ij}=\lambda^{(\epsilon(j)-\epsilon(i))/2}\overline{\sigma}_{-j,-i}~ \forall i,j\in\{1,...,-1\}$  and
\item $|\sigma_{*j}|\in\Lambda ~\forall j\in\{1,...,-1\}$. ($|\sigma_{*j}|=\sum\limits_{i=1}^n\bar\sigma_{ij}\sigma_{-i,j}$ is defined just before $3.4$.)
\end{enumerate}
}
\Proof{See \cite{bak-vavilov}, p.167.}
\Lemma{Let $\sigma\in U_{2n}(R,\Lambda)$, $x \in R^*$ and $k\in \{1,\dots,-1\}$. Then the statements below are true where $f_l:=e_l^t$ for any $l\in\{1,\dots,-1\}$.
\begin{enumerate}[(1)]
\item If the $k$-th column of $\sigma$ equals $x e_k$ then the $(-k)$-th row of $\sigma$ equals $\widebar{x^{-1}}f_{-k}$.
\item If the $k$-th row of $\sigma$ equals $x f_k$ then the $(-k)$-th column of $\sigma$ equals $\widebar{x^{-1}} e_{-k}$.
\end{enumerate}
}
\Proof{~\\
\vspace{-0.6cm}
\begin{enumerate}[(1)]
\item Since $\sigma^{-1}\sigma=e$ it follows that \[(\sigma^{-1}\sigma)_{ij}=\sum\limits_{l=1}^{-1}\sigma'_{il}\sigma_{lj}=\begin{cases} 1, & \text{if $i=j$,} 
\\
0, &\text{otherwise.}
\end{cases}\tag{3.1}\]
This implies that $1=\sum\limits_{l=1}^{-1}\sigma'_{kl}\sigma_{lk}=\sigma'_{kk}\sigma_{kk}=\sigma'_{kk}x$. Thus $\sigma'_{kk}=x^{-1}$. By Lemma 3.8, it follows that $\sigma_{-k,-k}=\widebar{x^{-1}}$. On the other hand $(3.1)$ implies that $0=\sum\limits_{l=1}^{-1}\sigma'_{il}\sigma_{lk}=\sigma'_{ik}\sigma_{kk}=\sigma'_{ik}x~\forall i\in\{1,\dots,-1\}\setminus\{k\}$. It follows that $\sigma'_{ik}=0~\forall i\in\{1,\dots,-1\}\setminus\{k\}$ and hence, by Lemma 3.8, $\sigma_{-k,-i}=0~\forall i\in\{1,\dots,-1\}\setminus\{k\}$, i.e. $\sigma_{-k,i}=0~\forall i\in\{1,\dots,-1\}\setminus\{-k\}$.   
\item Since $\sigma\sigma^{-1}=e$ it follows that \[(\sigma\sigma^{-1})_{ij}=\sum\limits_{l=1}^{-1}\sigma_{il}\sigma'_{lj}=\begin{cases} 1, & \text{if $i=j$,} 
\\
0, &\text{otherwise.}
\end{cases}\tag{3.2}\]
This implies $1=\sum\limits_{l=1}^{-1}\sigma_{kl}\sigma'_{lk}=\sigma_{kk}\sigma'_{kk}=x\sigma'_{kk}$. Thus $\sigma'_{kk}=x^{-1}$. By Lemma 3.8, it follows that $\sigma_{-k,-k}=\widebar{x^{-1}}$. On the other hand $(3.2)$ implies that $0=\sum\limits_{l=1}^{-1}\sigma_{kl}\sigma'_{lj}=\sigma_{kk}\sigma'_{kj}=x\sigma'_{kj}~\forall j\in\{1,\dots,-1\}\setminus\{k\}$. It follows that $\sigma'_{kj}=0~\forall j\in\{1,\dots,-1\}\setminus\{k\}$ and hence, by Lemma 3.8, $\sigma_{-j,-k}=0~\forall j\in\{1,\dots,-1\}\setminus\{k\}$, i.e. $\sigma_{j,-k}=0~\forall j\in\{1,\dots,-1\}\setminus\{-k\}$.\hfill$\qed$ 
\end{enumerate}}
\Definition{If $i,j\in\Omega$ such that $i\neq \pm j$ and $\xi\in R$, then the matrix
\[T_{ij}(\xi):=e+\xi e^{ij}-\lambda^{(\epsilon(j)-\epsilon(i))/2}\overline\xi e^{-j,-i}\in U_{2n}(R,\Lambda)\]
is called an {\it elementary short root element}. If $i\in \Omega$ and $\alpha\in\lambda^{-(\epsilon(i)+1)/2}\Lambda$, then the matrix \[T_{i,-i}(\alpha):=e+\alpha e^{i,-i}\in U_{2n}(R,\Lambda)\] is called an {\it elementary long root element}. If $\sigma\in U_{2n}(R,\Lambda)$ is an elementary short root element or an elementary long root element, it is called an {\it elementary unitary matrix}. The subgroup of $U_{2n}(R,\Lambda)$ generated by all elementary unitary matrices is called the {\it elementary unitary group} and is denoted by $EU_{2n}(R,\Lambda)$. Let $T_{ij}(\xi)$ be an elementary unitary matrix. If $i\neq -j\wedge\xi\in I$ or $i=-j \wedge \xi\in\lambda^{-(\epsilon(i)+1)/2}\Gamma$, then $T_{ij}(\xi)$ is called {\it elementary of level} $(I,\Gamma)$ or $(I,\Gamma)$-{\it elementary}. The subgroup of 
$U_{2n}(R,\Lambda)$ generated by all $(I,\Gamma)$-elementary matrices is called the {\it preelementary subgroup of level $(I,\Gamma)$} and is denoted by $EU_{2n}(I,\Gamma)$. Its normal closure in $EU_{2n}(R,\Lambda)$ is called the {\it elementary subgroup of level $(I,\Gamma)$} and is denoted by $EU_{2n}((R,\Lambda),(I,\Gamma))$. 
}
\Definition{Let $i,j\in\{1,\dots,-1\}$ such that $i\neq\pm j$. Define $P_{ij}:=e+e^{ij}-e^{ji}+\lambda^{(\epsilon(i)-\epsilon(j))/2}e^{-i,-j}-\lambda^{(\epsilon(j)-\epsilon(i))/2}e^{-j,-i}-e^{ii}-e^{jj}-e^{-i,-i}-e^{-j,-j}=T_{ij}(1)$ $T_{ji}(-1)T_{ij}(1)\in EU_{2n}(R,\Lambda)$. 
\Lemma{The relations \\
\begin{align*}T_{ij}(\xi)&=T_{-j,-i}(-\lambda^{(\epsilon(j)-\epsilon(i))/2}\overline{\xi})\tag{R1},\\
T_{ij}(\xi)T_{ij}(\zeta)&=T_{ij}(\xi+\zeta)\tag{R2},\\
[T_{ij}(\xi),T_{hk}(\zeta)]&=e,\tag{R3}\\
[T_{ij}(\xi),T_{jh}(\zeta)]&=T_{ih}(\xi\zeta),\tag{R4}\\
[T_{ij}(\xi),T_{j,-i}(\zeta)]&=T_{i,-i}(\xi\zeta-\lambda^{-\epsilon(i)}\bar\zeta\bar\xi)\text{ and}\tag{R5}\\
[T_{i,-i}(\alpha),T_{-i,j}(\xi)]&=T_{ij}(\alpha\xi)T_{-j,j}(-\lambda^{(\epsilon(j)-\epsilon(-i))/2)}\bar{\xi}\alpha\xi)\tag{R6}\end{align*}
hold where $h\neq j,-i$ and $k\neq i,-j$ in \textnormal{(R3)}, $i,h\neq\pm j$ and $i\neq \pm h$ in \textnormal{(R4)} and $i\neq\pm j$ in \textnormal{(R5)} and \textnormal{(R6)}.
}
\Proof{Straightforward calculation.}
\Definition{The group consisting of all $\sigma\in U_{2n}(R,\Lambda)$ such that $\sigma\equiv e(mod~I)$ and $\mathbbm{f}(\sigma u,\sigma u)\in \mathbbm{f}(u,u)+\Gamma~\forall u\in V$ is called the {\it principal congruence subgroup of level $(I,\Gamma)$} and is denoted by $U_{2n}((R,\Lambda),(I,\Gamma))$. By a theorem of Bak \cite{bak_2}, 4.1.4, cf. \cite{bak-vavilov}, 4.4, it is a normal subgroup of $U_{2n}(R,\Lambda)$.
}
\Lemma{Let $\sigma=\begin{pmatrix} a & b \\ c & d \end{pmatrix}\in U_{2n}(R,\Lambda)$, where $a,b,c,d\in M_n(R)$. Then $\sigma\in U_{2n}((R,\Lambda),(I,\Gamma))$ if and only if
\begin{enumerate}[(1)]
\item $\sigma\equiv e(mod~I)$ and
\item $|\sigma_{*j}|\in\Gamma ~\forall j\in\{1,...,-1\}$. ($|\sigma_{*j}|=\sum\limits_{i=1}^n\bar\sigma_{ij}\sigma_{-i,j}$ is defined just before $3.4$.)
\end{enumerate}
}
\Proof{See \cite{bak-vavilov}, p.174.}
\Definition{The preimage of the center of $U_{2n}(R,\Lambda)/U_{2n}((R,\Lambda),(I,\Gamma))$ under the canonical homomorphism $U_{2n}(R,\Lambda)\rightarrow U_{2n}(R,\Lambda)/U_{2n}((R,\Lambda),(I,\Gamma))$ is called the {\it full congruence subgroup of level $(I,\Gamma)$} and is denoted by $CU_{2n}((R,\Lambda),(I,\Gamma))$.
}
\Remark{Obviously $U_{2n}((R,\Lambda),(I,\Gamma))\subseteq CU_{2n}((R,\Lambda),(I,\Gamma))$ and $CU_{2n}((R,\Lambda),$ $(I,\Gamma))$ is a normal subgroup of $U_{2n}(R,\Lambda)$.}
\Lemma{If $n\geq 3$ and $R$ is almost commutative (i.e. module finite over its center), then the equalities
\begin{align*}&[CU_{2n}((R,\Lambda),(I,\Gamma)),EU_{2n}(R,\Lambda)]\\
=&[EU_{2n}((R,\Lambda),(I,\Gamma)),EU_{2n}(R,\Lambda)]\\
=&EU_{2n}((R,\Lambda),(I,\Gamma))\end{align*}
hold.
\Proof{See \cite{bak-vavilov}, Theorem 1.1 and Lemma 5.2.\hfill$\qed$}
\section[Sandwich classification for hyperbolic unitary groups]{Sandwich classification for hyperbolic unitary\newline groups}
In this section we prove our main results.\\

We begin by fixing notation which will be used for most of the section, up through the proof of Theorem 4.11. This theorem proves our main result over a special kind of almost commutative form ring, which will be described immediately below. The general result over quasi-finite form rings will be deduced from this result in Theorem 4.14. After fixing notation below, we shall explain the idea of the proof of Theroem 4.11 and how the rest of the section is organized.\\

Until the end of the proof of Theorem 4.11 let $n\geq 3$, $(R,\Lambda)$ a form ring and $C$  the subring of $R$ consisting of all finite sums of elements of the form $c\bar c$ and $-c\bar c$ where $c$ ranges over some subring $C'\subseteq center(R)$. Further assume that $R$ is a Noetherian $C$-module. For any form ideal $(I,\Gamma)$ of $(R,\Lambda)$ and maximal ideal $m$ of $C$ set $S_m:=C\setminus m$, $R_m:=S_m^{-1}R$, $\Lambda_m:=S_m^{-1}\Lambda$, $I_m:=S_m^{-1}I$ and $\Gamma_m:=S_m^{-1}\Gamma$. Let \[\phi_m: U_{2n}(R,\Lambda)/U_{2n}((R,\Lambda),(I,\Gamma))\rightarrow U_{2n}(R_m,\Lambda_m)/U_{2n}((R_m,\Lambda_m),(I_m,\Gamma_m))\] be the homomorphism induced by $F_m$ where \[F_m:U_{2n}(R,\Lambda)\rightarrow U_{2n}(R_m,\Lambda_m)\] is the homomorphism induced by the localisation homomorphism \[f_m:R\rightarrow R_m.\]Let \[\psi: U_{2n}(R,\Lambda)\rightarrow U_{2n}(R,\Lambda)/U_{2n}((R,\Lambda),(I,\Gamma))\] and \[\psi_m: U_{2n}(R_m,\Lambda_m)\rightarrow U_{2n}(R_m,\Lambda_m)/U_{2n}((R_m,\Lambda_m),
(I_m,\Gamma_m))\]be the canonical homomorphisms. Further set $\lambda_m:=f_m(\lambda)$. Note that the diagram
\xymatrixcolsep{3pc}
\xymatrixrowsep{4pc}
\[\xymatrix{
U_{2n}(R, \Lambda) \ar[d]^{F_m} \ar[r]^/-2.1em/{\psi} &U_{2n}(R,\Lambda)/U_{2n}((R,\Lambda),(I,\Gamma))\ar[d]^{\phi_m}\\
U_{2n}(R_m, \Lambda_m) \ar[r]^/-2.1em/{\psi_m} &U_{2n}(R_m,\Lambda_m)/U_{2n}((R_m,\Lambda_m),(I_m,\Gamma_m))}\]
is commutative for any form ideal $(I,\Gamma)$ of $(R,\Lambda)$ and maximal ideal $m$ of $C$.\\

A road map of the proof of Theorem 4.11 is as follows. It is easy to show that there is a largest form ideal $(I,\Gamma)$ of $(R,\Lambda)$ such that $EU_{2n}((R,\Lambda), (I,\Gamma)) \subseteq H$. To show that $H \subseteq CU_{2n}((R,\Lambda), (I,\Gamma))$ we proceed by contradiction. We show that if $H \not\subseteq  CU_{2n}((R,\Lambda), (I,\Gamma))$, then there is a form ideal $(I',\Gamma')$ which properly contains $(I,\Gamma)$ such that $EU_{2n}((R,\Lambda), (I',\Gamma')) \subseteq H\cdot U_{2n}((R,\Lambda),(I,\Gamma))$. Then it follows routinely from the standard mixed commutator formulas of \cite{bak-vavilov}, that $EU_{2n}((R,\Lambda), (I',\Gamma')) \subseteq H$, which contradicts the maximality of $(I,\Gamma)$. This is the approach pioneered by W. Klingenberg and H. Bass for the general linear group $GL_n$ over respectively semilocal rings and rings $R$ satisfying a stability condition and can also be applied over almost commutative rings. Over the last rings, the crucial step in the road map is embedding $GL_n(R)/GL_n(R,I)$ in the canonical way into $GL_n(R/I)$ and then showing that if $H'=$ image of $H$ in $GL_n(R/I)$ is noncentral in $GL_n(R/I)$, then $H'$ contains a relative elementary group $E_n(R/I,I'/I)$ for some ideal $I'$ of $R$, which properly contains $I$. This result is established by a localization argument over the maximal ideals of $center(R/I)$. In the case of unitary groups, we would like to apply the same approach, but unfortunately this doesn't work, because for an arbitrary form ideal $(I,\Gamma)$, there is no unitary group in the sense of the current paper into which we can canonically embed $U_{2n}(R,\Lambda)/U_{2n}((R,\Lambda),(I,\Gamma))$. But, there is an unpublished theory of "model unitary groups" of A. Bak, which has the  property that the quotient of any model unitary group by one of its congruence subgroups is again a model unitary group and which has an effective localization procedure. Of course, any unitary group $U_{2n}(R,\Lambda)$ is a model unitary group, but not conversely. In this section, we shall adapt the methods of model unitary groups to the circumstances at hand and show directly that the quotient group $U_{2n}(R,\Lambda)/U_{2n}((R,\Lambda),(I,\Gamma))$ has a good localization theory. The essential contributions here are made in Lemmas 4.2 and 4.4. Lemma 4.8 is an application of this localization theory and provides the key result needed to show that if $\psi(H) \not\subseteq center(U_{2n}(R,\Lambda)/U_{2n}((R,\Lambda),(I,\Gamma)))$ then there is a form ideal $(I',\Gamma')$ which properly contains $(I,\Gamma)$ such that $EU_{2n}((R,\Lambda)),(I',\Gamma')) \subseteq H\cdot U_{2n}((R,\Lambda),(I,\Gamma))$, i.e $\psi(EU_{2n}((R,\Lambda)),(I',\Gamma')))\subseteq \psi(H)$. The proof of Theorem 4.11 can be then concluded as above by applying the standard mixed commutator formulas.
\Lemma{Let $(I,\Gamma)$ be a form ideal of $(R,\Lambda)$ and $g'\in U_{2n}(R,\Lambda)/U_{2n}((R,\Lambda),(I,$ $\Gamma))$ be noncentral. Then there is a maximal ideal $m$ of $C$ such that $I\cap C\subseteq m$ and $\phi_m(g')$ is noncentral.}
\Proof{Since $g'\in U_{2n}(R,\Lambda)/U_{2n}((R,\Lambda),(I,\Gamma))$ is noncentral, there is an $h'\in U_{2n}(R,$ $\Lambda)/U_{2n}((R,\Lambda),(I,\Gamma))$ such that $g'h'\neq h'g'$. Let $g,h\in U_{2n}(R,\Lambda)$ such that $g'=gU_{2n}((R,\Lambda),(I,
\Gamma))$ and $h'=hU_{2n}((R,\Lambda),(I,\Gamma))$. Set $\sigma:=[g^{-1},h^{-1}]$. Clearly $g'h'\neq h'g'$ implies $\sigma\not\in U_{2n}((R,\Lambda),$ $(I,\Gamma))$. Hence either $\sigma_{ij}\not\in I$ for some $i,j\in\{1,\dots,-1\}$ such that $i\neq j$, or $\sigma_{ii}-1\not\in I$ for some $i\in\{1,\dots,-1\}$ or $x_j:=|\sigma_{*j}|\not \in\Gamma$ for some $j\in\{1,\dots,-1\}$.\\
\\
\underline{case 1} Assume that $\sigma_{ij}\not\in I$ for some $i,j\in\{1,\dots,-1\}$ such that $i\neq j$. Set $Y:=\{c\in C|c\sigma_{ij}\in I\}$. Since $\sigma_{ij}\not\in I$, $Y$ is a proper ideal of $C$. Hence it is contained in a maximal ideal $m$ of $C$. Clearly $I\cap C\subseteq Y\subseteq m$ and hence $S_m\cap Y=\emptyset$. We show now that $\phi_m(g')$ does not commute with $\phi_m(h')$, i.e. $F_m(\sigma)\not\in U_{2n}((R_{m},\Lambda_{m}),(I_{m},\Gamma_{m}))$. Obviously $(F_m(\sigma))_{ij}=f_m(\sigma_{ij})$. Assume $(F_m(\sigma))_{ij}\in I_{m}$. Then
\begin{align*}
&\exists x\in I, s\in S_m: \frac {\sigma_{ij}}{1}=\frac{x}{s}\\
\Rightarrow~&\exists x\in I, s,t\in S_m: t(\sigma_{ij}s-x)=0\\
\Rightarrow~&\exists x\in I, s,t\in S_m: st\sigma_{ij}=tx\in I\\
\Rightarrow~&\exists u\in S_m: u\sigma_{ij}\in I.
\end{align*}
But this contradicts $S_m\cap Y=\emptyset$. Hence $(F_m(\sigma))_{ij}\not\in I_m$ and thus $\phi_m(g')$ is noncentral.\\
\\
\underline{case 2} Assume that $\sigma_{ii}-1\not\in I$ for some $i\in\{1,\dots,-1\}$. Set $Y:=\{c\in C|c(\sigma_{ii}-1)\in I\}$. Since $\sigma_{ii}-1\not\in I$, $Y$ is a proper ideal of $C$. Hence it is contained in a maximal ideal $m$ of $C$. Clearly $I\cap C\subseteq Y\subseteq m$ and hence $S_m\cap Y=\emptyset$. We show now that $\phi_m(g')$ does not commute with $\phi_m(h')$, i.e. $F_m(\sigma)\not\in U_{2n}((R_{m},\Lambda_{m}),(I_{m},\Gamma_{m}))$. Obviously $(F_m(\sigma))_{ii}-1=f_m(\sigma_{ii})-1=f_m(\sigma_{ii}-1)$. Assume $(F_m(\sigma))_{ii}-1\in I_{m}$. Then
\begin{align*}
&\exists x\in I, s\in S_m: \frac {\sigma_{ii}-1}{1}=\frac{x}{s}\\
\Rightarrow~&\exists x\in I, s,t\in S_m: t((\sigma_{ii}-1)s-x)=0\\
\Rightarrow~&\exists x\in I, s,t\in S_m: st(\sigma_{ii}-1)=tx\in I\\
\Rightarrow~&\exists u\in S_m: u(\sigma_{ii}-1)\in I.
\end{align*}
But this contradicts $S_m\cap Y=\emptyset$. Hence $(F_m(\sigma))_{ii}-1\not\in I_{m}$ and thus $\phi_m(g')$ is noncentral.\\
\\
\underline{case 3} Assume that $x_j=|\sigma_{*j}|\not \in\Gamma$ for some $j\in\{1,\dots,-1\}$. Set $Y:=\{c\in C|cx_j\in \Gamma\}$. Since $x_j\not\in \Gamma$, $Y$ is a proper ideal of $C$. Hence it is contained in a maximal ideal $m$ of $C$. Since $x_j\in \Lambda$ and $y^2\Lambda\subseteq\Gamma_{min}\subseteq\Gamma$ for any $y\in I\cap C$, $(I\cap C)^2\subseteq Y\subseteq m$. This implies $S_m\cap Y=\emptyset$ and $I\cap C\subseteq m$, since $m$ is prime. We show now that $\phi_m(g')$ does not commute with $\phi_m(h')$, i.e. $F_m(\sigma)\not\in U_{2n}((R_{m},\Lambda_{m}),(I_{m},\Gamma_{m}))$. Obviously $|(F_m(\sigma))_{*j}|=f_m(x_j)$. Assume $f_m(x_j)\in \Gamma_{m}$. Then
\begin{align*}
&\exists y\in \Gamma, s\in S_m: \frac {x_j}{1}=\frac{y}{s}\\
\Rightarrow~&\exists y\in \Gamma, s,t\in S_m: t(x_js-y)=0\\
\Rightarrow~&\exists y\in \Gamma, s,t\in S: stx_j=ty\in \Gamma\\
\Rightarrow~&\exists u\in S_m: ux_j\in \Gamma.
\end{align*}
But this contradicts $S_m\cap Y=\emptyset$. Hence $|(F_m(\sigma))_{*j}|=f_m(x_j)\not\in \Gamma_{m}$ and thus $\phi_m(g')$ is noncentral.
\hbox{}\hfill$\qed$
}
\Lemma{Let $(I,\Gamma)$ be a form ideal of $(R, \Lambda)$ and $m$ a maximal ideal of $C$. Then there is an $s_0\in S_m$ with the properties 
\begin{enumerate}[(1)]
\item if $x \in s_0R$ and $\exists t\in S_m: tx\in I$, then $x\in I$ and
\item if $x \in s_0R$ and $\exists t\in S_m: tx\in \Gamma$, then $x\in \Gamma$.
\end{enumerate}
It follows that $\phi_m$ is injective on $\psi(U_{2n}((R,\Lambda),(s_0R,s_0\Lambda)))$.
}
\Proof{For any $s\in S_m$ set $Y(s):=\{x\in R|sx\in I\}$. Then for any $s\in S_m$, $Y(s)$ is a $C$-submodule of $R$. Since $R$ is Noetherian $C$-module, the set $\{Y(s)|s\in S_m\}$ has a maximal element $Y(s_1)$. Clearly all elements $x\in s_1R$ have the property that $tx\in I$ for some $t\in S_m$ implies $x\in I$. For any $s\in S$ set $Z(s):=\{x\in R|sx\in \Gamma\}$. Then for any $s\in S_m$, $Z(s)$ is a $C$-submodule of $R$. Since $R$ is a Noetherian $C$-module, the set $\{Z(s)|s\in S\}$ has a maximal element $Z(s_2)$. Clearly all elements $x\in s_2R$ have the property that $tx\in \Gamma$ for some $t\in S_m$ implies $x\in \Gamma$. Set $s_0:=s_1s_2$. Since $s_0R=s_1s_2R\subseteq s_1R\cap s_2R$, $s_0$ has the properties $(1)$ and $(2)$ above. We will show now that $\phi_m$ is injective on $\psi(U_{2n}((R,\Lambda),(s_0R,s_0\Lambda)))$.  Let $g'_1,g'_2\in \psi(U_{2n}((R,\Lambda),(s_0R,s_0\Lambda)))$ such that $\phi_m(g'_1)=\phi_m(g_2')$. Since $g'_1,g'_2\in \psi(U_{2n}((R,\Lambda),(s_0R,s_0\Lambda)
))$, there are $g_1,g_2\in 
U_{2n}((R,\Lambda),(s_0R,s_0\Lambda))$ such that $\psi(g_1)=g'_1$ and $\psi(g_2)=g'_2$. Set $h:=(g_1)^{-1}g_2\in U_{2n}((R,\Lambda),(s_0R,s_0\Lambda))$. Clearly $\phi_m(g'_1)=\phi_m(g_2')$ is equivalent to $F_m(h)\in U_{2n}((R_m,\Lambda_m),$ $(I_m,\Gamma_m))$, i.e.
\begin{enumerate}[(a)]
\item $F_m(h)\equiv e (mod~I_m)$ and
\item $f_m(|h_{*j}|)\in\Gamma_m~\forall j\in\{1,\dots,-1\}$.
\end{enumerate}
We want to show that $g'_1=g'_2$ which is equivalent to $h\in U_{2n}((R,\Lambda),(I,\Gamma))$, i.e.
\begin{enumerate}[(a')]
\item $h\equiv e (mod~I)$ and
\item $x_j:=|h_{*j}|\in\Gamma~\forall j\in\{1,\dots,-1\}$.
\end{enumerate}
First we show (a'). Let $i,j\in\{1,\dots,-1\}$ such that $i\neq j$. Since (a) holds, $f_m(h_{ij})\in I_m$. Hence 
\begin{align*}
&\exists x\in I, s\in S_m: \frac {h_{ij}}{1}=\frac{x}{s}\\
\Rightarrow~&\exists x\in I, s,t\in S_m: t(h_{ij}s-x)=0\\
\Rightarrow~&\exists x\in I, s,t\in S_m: sth_{ij}=tx\in I\\
\Rightarrow~&\exists u\in S_m: uh_{ij}\in I.\tag{4.1.1}
\end{align*}
Since $h\in U_{2n}((R,\Lambda),(s_0R,s_0\Lambda))$, $h_{ij}\in s_0R$. It follows from $(4.1.1)$ that $h_{ij}\in I$ since $s_0$ has property $(1)$. Analogously one can show that $h_{ii}-1\in I$ for all $i\in\{1,\dots,-1\}$. Hence $h\equiv e (mod~I)$. Now we show (b'). Let $j\in\{1,\dots,-1\}$. Since (b) holds, $f_m(x_j)\in\Gamma_m$. Hence
\begin{align*}
&\exists y\in \Gamma, s\in S_m: \frac {x_j}{1}=\frac{y}{s}\\
\Rightarrow~&\exists y\in \Gamma, s,t\in S_m: t(x_js-y)=0\\
\Rightarrow~&\exists y\in \Gamma, s,t\in S_m: stx_j=ty\in \Gamma\\
\Rightarrow~&\exists u\in S_m: ux_j\in \Gamma.\tag{4.1.2}
\end{align*}
Since $h\in U_{2n}((R,\Lambda),(s_0R,s_0\Lambda))$, $x_j\in s_0\Lambda$. It follows from $(4.1.2)$ that $x_{j}\in \Gamma$ since $s_0$ has property $(2)$. Hence $g'_1=g'_2$ and thus $\phi_m$ is injective on $\psi(U_{2n}((R,\Lambda),(s_0R,$ $s_0\Lambda)))$. \hbox{}\hfill$\qed$}
\Definition{Let $G$ denote a group and $A$ a set of subgroups of $G$ such that
\begin{enumerate}[(1)]
\item for any $U,V\in A$ there is a $W\in A$ such that $W\subseteq U\cap V$ and
\item for any $g\in G$ and $U\in A$ there is a $V\in A$ such that $^gV\subseteq U$.
\end{enumerate} 
Then $A$ is called a {\it base of open subgroups of $1\in G$}. A pair $(A,B)$ is called a
{\it supplemented base} for $G$ if $A$ and $B$ and are sets of nontrivial
subgroups of $G$ such that $A$ is a base
of open subgroups of $1\in G$, each member of $B$ is contained in some member of $A$,
and if $U \in A$ and $V \in B$ then $U\cap V$ contains a member of $B$.}\\

In the lemma below we use the following conventions. Let $x\in R$. Then $RxR$ denotes the {\it involution invariant ideal generated by $x$}, i.e. the ideal of $R$ 
generated by $\{x,\bar x\}$. Now let $(I,\Gamma)$ be a form ideal of $(R,\Lambda)$ and assume that $x\in R\setminus I$ or $x\in \Gamma^I_{max}\setminus\Gamma$. Set $\Gamma(x):=\Gamma^{RxR}_{min}$ if $x\in R\setminus I$ and $\Gamma(x):=\Gamma^{RxR}_{min}+\langle\{yx\bar y|y\in R\}\rangle$ if $x\in \Gamma^I_{max}\setminus\Gamma$. $\Gamma(x)$ is called the {\it relative form parameter defined by $x$ and $(I,\Gamma)$}. One checks easily that $(RxR,\Gamma(x))$ is a form ideal of $(R,\Lambda)$ which is not contained in $(I,\Gamma)$, i.e. $RxR\not\subseteq I$ or $\Gamma(x)\not\subseteq\Gamma$. It is called the {\it form ideal defined by $x$ and $(I,\Gamma)$}.
\Lemma{Let $(I,\Gamma)$ be a form ideal of $(R,\Lambda)$, $m$ a maximal ideal of $C$ and $s_0\in S_m$ as in the previous lemma. Set $A:=\{EU_{2n}(ss_0R,ss_0\Lambda)|s\in S_m\}$ and $B:=\{EU_{2n}(Rxs_0R,$ $\Gamma(xs_0))|(x\in R, xs_0\in R\setminus I) \lor (x\in \Lambda, xs_0\in I\setminus\Gamma)\}$. Then $(A,B)$ is a supplemented base for $EU_{2n}(R,\Lambda)$ and $F_m(A,B):=(F_m(A),F_m(B))$ is a supplemented base for $EU_{2n}(R_m,\Lambda_m)$.
}
\Proof{First we show $(A,B)$ is a supplemented base for $EU_{2n}(R,\Lambda)$. Clearly $A$ and $B$ are sets of nontrivial subgroups of $EU_{2n}(R,\Lambda)$. We show now that $A$ is a base of open subgroups of $1\in EU_{2n}(R,\Lambda)$. Therefore we must show that $A$ satisfies the conditions $(1)$ and $(2)$ in Definition 4.3.
\begin{enumerate}[(1)]
\item Let $U=EU_{2n}(ss_0R,ss_0\Lambda),V=EU_{2n}(ts_0R,ts_0\Lambda)\in A$. Set $W:=EU_{2n}(sts_0R,$ $sts_0\Lambda)\in A$. Then clearly $W\subseteq U\cap V$.
\item Let $g\in EU_{2n}(R,\Lambda)$ and $U=EU_{2n}(ss_0R,ss_0\Lambda)\in A$. There is a $K\in\mathbb{N}$ such that $g$ is the product of $K$ elementary unitary matrices. Set $V:=EU_{2n}((ss_0)^{2\cdot 4^K+}$ $^{4^{K-1}+\dots+4}R,(ss_0)^{2\cdot 4^K+4^{K-1}+\dots+4}\Lambda)\in A$. Then $^gV\subseteq U$ (see Lemma 4.1 in \cite{hazrat}).
\end{enumerate}
Hence $A$ is a base of open subgroups of $1\in EU_{2n}(R,\Lambda)$. Let $EU_{2n}(Rxs_0R,\Gamma(xs_0))\in B$. Then $EU_{2n}(Rxs_0R,\Gamma(xs_0))\subseteq EU_{2n}(s_0R,s_0\Lambda)\in A$. It remains to show that if $U\in A$ and $V\in B$ then $U\cap V$ contains a member of $B$. Let $U=EU_{2n}(ss_0R,ss_0\Lambda)\in A$ and $V=EU_{2n}(Rxs_0R,\Gamma(xs_0))\in B$. Set $W:=EU_{2n}(Rxss_0R,\Gamma(xss_0))$. If $xs_0\not\in I$, then $xss_0\not\in I$ and if $xs_0\not\in\Gamma$, then $xss_0\not \in\Gamma$ (by the definition of $s_0$, see the previous lemma). Hence $W\in B$. Obviously $W\in U\cap V$. Hence $(A,B)$ is a supplemented base for $EU_{2n}(R,\Lambda)$.\\

Now we show $F_m(A,B)$ is a supplemented base for $EU_{2n}(R_m,\Lambda_m)$. Clearly $F_m(A)$ and $F_m(B)$ are sets of nontrivial subgroups of $EU_{2n}(R_m,\Lambda_m)$. We show now that $F_m(A)$ is a base of open subgroups of $1\in EU_{2n}(R_m,\Lambda_m)$. Therefore we must show that $F_m(A)$ satisfies the conditions $(1)$ and $(2)$ in Definition 4.3.
\begin{enumerate}[(1)]
\item Let $U=F_m(EU_{2n}(ss_0R,ss_0\Lambda)),V=F_m(EU_{2n}(ts_0R,ts_0\Lambda))\in F_m(A)$. Set $W:=F_m(EU_{2n}(sts_0R,sts_0\Lambda))\in F_m(A)$. Then clearly $W\subseteq U\cap V$.
\item Let $g\in EU_{2n}(R_m,\Lambda_m)$ and $U=F_m(EU_{2n}(ts_0R,ts_0\Lambda))\in F_m(A)$. There are a $K\in\mathbb{N}$ and elementary unitary matrices $\tau_1=T_{i_1j_1}(\frac{x_1}{s_1}),\dots,\tau_K=T_{i_kj_k}(\frac{x_K}{s_K})\in EU_{2n}(R_m,\Lambda_m)$ such that $g=\tau_1\dots\tau_K$. Set $s:=s_1\dots s_K$ and  $V:=F_m(EU_{2n}($ $(sts_0)^{2\cdot 4^K+4^{K-1}+\dots+4}R,(sts_0)^{2\cdot 4^K+4^{K-1}+\dots+4}\Lambda))\in F_m(A)$. Then $^gV\subseteq U$ (see Lemma 4.1 in \cite{hazrat}).
\end{enumerate}
Hence $F_m(A)$ is a base of open subgroups of $1\in EU_{2n}(R_m,\Lambda_m)$. That each member of $F_m(B)$ is contained in some member of $F_m(A)$ follows from the fact that any member of $B$ is contained in a member of $A$. That given $U\in F_m(A)$ and $V\in F_m(B)$, $U\cap V$ contains a member of $F_m(B)$ follows from the fact that given $U\in A$ and $V\in B$, $U\cap V$ contains a member of $B$. Hence $F_m(A,B)$ is a supplemented base for $EU_{2n}(R_m,\Lambda_m)$. \hfill $\qed$
}\\

The lemmas 4.5, 4.6 and 4.7 will be used in the proof of Lemma 4.8.
\Lemma{Let $(I,\Gamma)$ be a form ideal of $(R,\Lambda)$ and $S\subseteq C$ a multiplicative subset. Let $T_{ij}(x)\in EU_{2n}(R_m,\Lambda_m)$ be an elementary short or long root element, $\sigma\in U_{2n}(R_m,\Lambda_m)$ and $s\in S$. Then $[\sigma,T_{ij}(x)]\in U_{2n}((R_m,\Lambda_m),(I_m,\Lambda_m\cap I_m))$ if and only if $[\sigma,T_{ij}(f_m(s)x)] \in U_{2n}((R_m,\Lambda_m),(I_m,\Lambda_m\cap I_m))$.
}
\Proof{Straightforward computation.\hfill$\qed$}
\Lemma{Let $\sigma\in U_{2n}(R,\Lambda)$. Further let $i,j\in \{1,\dots,-1\}$ such that $i\neq \pm j$, $x\in R$ and $y\in \lambda^{-(\epsilon(i)+1)/2}\Lambda$. Set $\tau:=[\sigma,T_{ij}(x)]$ and $\rho:=[\sigma,T_{i,-i}(y)]$. Then 
\[|\tau_{*k}|=\bar\sigma'_{jk}\bar x|\sigma_{*i}|x\sigma'_{jk}+\bar\sigma'_{-i,k} x |\sigma_{*,-j}|\bar x \sigma'_{-i,k}+a_k-\lambda\bar a_k,\]
if $k\neq j,-i$,
\[|\tau_{*j}|=\bar\sigma'_{jj}\bar x|\sigma_{*i}|x\sigma'_{jj}+\bar\sigma'_{-i,j} x |\sigma_{*,-j}|\bar x \sigma'_{-i,j}+\bar x|\tau_{*i}|x+a_j-\lambda\bar a_j\]
and
\[|\tau_{*,-i}|=\bar\sigma'_{j,-i}\bar x|\sigma_{*i}|x\sigma'_{j,-i}+\bar\sigma'_{-i,-i} x |\sigma_{*,-j}|\bar x \sigma'_{-i,-i}+x|\tau_{*,-j}|\bar x+a_{-i}-\lambda\bar a_{-i}\]
where for each $k\in\{1,\dots,-1\}$, $a_k$ lies in the ideal $J(\sigma)$ generated by the nondiagonal entries of $\sigma$ and $\sigma^{-1}$.
It follows that if $I$ is an involution invariant ideal and $\sigma\in U_{2n}((R,\Lambda),(I,I\cap\Lambda))$, then $|\tau_{*k}|\in \Gamma^I_{min}~\forall k\neq j,-i$, $|\tau_{*j}|\equiv\bar x|\sigma_{*i}|x(\textnormal{mod } \Gamma^I_{min})$ and $|\tau_{*,-i}|\equiv x|\sigma_{*,-j}|\bar x(\textnormal{mod } \Gamma^I_{min})$. Further
\[|\rho_{*k}|=\bar\sigma'_{-i,k}\bar y|\sigma_{*i}|y\sigma'_{-i,k}+b_k,\]
if $k\neq -i$,
\[|\rho_{*,-i}|=\bar\sigma'_{-i,-i}\bar y|\sigma_{*i}|y\sigma'_{-i,-i}+\bar y|\rho_{*i}|y+b_{-i}+c-\lambda\bar c\]
where
\[
b_k=\begin{cases}
\sigma_{-k,i}y\sigma'_{-i,k},&\text{if }k\neq -i \text{ and }\epsilon(k)=1,\\
\overline{\sigma_{-k,i}y\sigma'_{-i,k}},&\text{if }k\neq -i \text{ and }\epsilon(k)=-1,\\
\overline{\sigma_{ii}y\sigma'_{-i,i}-y},&\text{if }k= -i \text{ and }\epsilon(i)=1,\\
\sigma_{ii}y\sigma'_{-i,-i}-y,&\text{if }k=-i \text{ and }\epsilon(i)=-1,\\
\end{cases}
\]
and $c\in J(\sigma)$. It follows that if $I$ is an involution invariant ideal and $\sigma\in U_{2n}((R,\Lambda),$ $(I,I\cap\Lambda))$, then $|\rho_{*k}|\in \Gamma^I_{min}~\forall k\neq -i$ and $|\rho_{*,-i}|\equiv \bar y|\sigma_{*i}|y(\textnormal{mod } \Gamma^I_{min})$.
}
\Proof{Straightforward computation.\hfill$\qed$}\\
\Lemma{Let $m$ be a maximal ideal of $C$ and $\sigma\in U_{2n}(R_m,\Lambda_m)$. Then there is an $\epsilon\in EU_{2n}(R_m,\Lambda_m)$ such that $(^{\epsilon}\sigma)_{11}$ is invertible.}
\Proof{By Lemma 1.4 in \cite{vaserstein-you} and Lemma 3.4 in \cite{bak-tang}, $R_m$ satisfies the $\Lambda$-stable range condition $\Lambda S_1$. Hence there is an $\epsilon_1=\begin{pmatrix}e^{n\times n}&0\\\gamma&e^{n\times n}\end{pmatrix}\in EU_{2n}(R_m,\Lambda_m)$, where $\gamma\in M_{n}(R_m)$, such that $(x_1,\dots,x_n)$ is right unimodular where $(x_1,\dots,x_{-1})$ is the first row of $^{\epsilon_1}\sigma$. Since $\Lambda S_1$ implies $S R_1$, there is a matrix
$\epsilon_2=\begin{pmatrix}\omega_1&0\\0&\omega_2\end{pmatrix}\in EU_{2n}(R_m,\Lambda_m)$, where $\omega_1$ and $\omega_2$ are lower triangular matrices in $M_{n}(R_m)$ with $1$'s on the diagonal, such that the entry of $(^{\epsilon_2\epsilon_1}\sigma)_{11}$ is right invertible. Since $R$ is a Noetherian $C$-module, $R$ is almost commutative. It follows that $R_m$ is almost commutative and hence $(^{\epsilon_2\epsilon_1}\sigma)_{11}$ is invertible (note that almost commutative rings are Dedekind-finite by Nakayama's Lemma). \hfill $\qed$ 
}\\

In the following lemma we will apply lemmas and corollaries in \cite{bak_2}, chapter IV, \S 3. We are allowed to do this since for any maximal ideal $m$ of $C$, $C'_m:=S_m^{-1}C'$ is semilocal by Lemma 1.4 in \cite{vaserstein-you} and hence the Bass-Serre-dimension of $C'_m$ is $0$. Since $R$ is module finite over $C'$, $R_m$ is module finite over $C'_m$ 
and hence $R_m$ is a finite $C'_m$-algebra. Further set $A':=\phi_m(\psi(A))$ and $B':=\phi_m(\psi(B))$ where $(A,B)$ is the supplemented base for $EU_{2n}(R,\Lambda)$ defined in Lemma 4.4. Clearly $(A',B')$ is a supplemented base for $\psi_m(EU_{2n}(R_m,\Lambda_m))$ since $F_m(A,B)$ is a supplemented base for $EU_{2n}(R_m,\Lambda_m)$ by Lemma 4.4 and $\phi_m\circ \psi=\psi_m\circ F_m$.
\Lemma{Let $(I,\Gamma)$ be a form ideal of $(R,\Lambda)$, $m$ a maximal ideal of $C$ and $h'\in U_{2n}(R_m,\Lambda_m)/U_{2n}((R_m,\Lambda_m),(I_m,\Gamma_m))$ a noncentral element. Then given any $U'\in A'$, there is a $k\in\mathbb{N}$ and elements $g'_0,\dots,g'_k\in\phi_m(\psi(EU_{2n}(R,\Lambda)))$, $\epsilon'_0,\dots,\epsilon'_k\in \psi_m(EU_{2n}(R_m,$  $\Lambda_m))$ and $l_1,\dots,l_k\in\{-1,1\}$ such that $g'_k$ is the nontrivial image of an elementary matrix in $EU_{2n}(R,\Lambda)$, 
\[ ^{\epsilon'_{k}}([^{\epsilon'_{k-1}}(\dots^{\epsilon'_2}([^{\epsilon'_1}([^{\epsilon'_0}h',g'_0]^{l_1}),g'_1]^{l_2})\dots),g'_{k-1}]^{l_{k}})=g'_k\]
and
\[^{d'_i}g'_i\in U'~\forall i\in \{0,\dots,k\}\]
where $d'_i=(\epsilon'_{i}\cdot\hdots\cdot\epsilon'_0)^{-1}\hspace{0.1cm}~(0\leq i\leq k)$. 
}
\Proof{Let $U'\in A'$. Then there is a $U\in F_m(A)$ such that $\psi_m(U)=U'$. Let $h\in U_{2n}(R_m,\Lambda_m)$ such that $h'=\psi_m(h)$. Since $h'$ is noncentral, $h\not\in CU_{2n}((R_m,\Lambda_m),(I_m,$ $\Gamma_m))$. The proof is divided into three parts, I, II and III. In Part I we assume that $h\not\in CU_{2n}((R_m,\Lambda_m),(I_m,I_m\cap\Lambda_m))$ and $n>3$. In Part II we assume that $h\not\in CU_{2n}((R_m,\Lambda_m),(I_m,I_m\cap\Lambda_m))$ and $n=3$. In Part III we assume that $h\in CU_{2n}((R_m,$ $\Lambda_m),(I_m,I_m\cap\Lambda_m))$. The groups $U_i\in F_m(A)~(0\leq i\leq k)$ appearing in the proof are chosen such that $^{d_i}U_i\subseteq U$ where $d_i=(\epsilon_{i}\cdot\hdots\cdot\epsilon_0)^{-1}$ (possible by Lemma 4.4). The elements $t_{i}\in S_m~(0\leq i\leq k)$ are chosen such that $U_i=F_m(EU_{2n}(t_is_0R,t_is_0\Lambda))$. Further we denote $f_m(t_is_0)$ by $s_i~(0\leq i\leq k)$.\\
\\
\underline{Part I} Assume that $h\not\in CU_{2n}((R_m,\Lambda_m),(I_m,I_m\cap\Lambda_m))$ and $n>3$.\\
By \cite{bak_2}, chapter IV, Corollary 3.10 (applied with $H=$$^{EU_{2n}(R_m,\Lambda_m)}\langle h\rangle$), there is an $\epsilon_0\in EU_{2n}(R_m,\Lambda_m)$ and an $x=\frac{a}{s}
\in R_m$ such that $[^{\epsilon_0} h,T_{1,-2}(x)]\not\in CU_{2n}((R_m,\Lambda_m),$ $(I_m,I_m\cap\Lambda_m))$. By \cite{bak_2}, chapter IV, Lemma 3.12, Part I, case 7 there is a matrix $\epsilon_1\in EU_{2n}(R_m,\Lambda_m)$ of the form \[\epsilon_1=\begin{pmatrix}X & Y \\ 0 & Z\end{pmatrix},\] where $X,Y,Z\in M_n(R_m)$, such that the first $n$ coordinates of $\epsilon_1(^{\epsilon_0} h)_{*1}$ equal $\begin{pmatrix} 1&0&\dots&0 \end{pmatrix}^{t}$ and the first $n$ coordinates of $\epsilon_1(^{\epsilon_0} h)_{*2}$ equal $\begin{pmatrix} 0&1&0&\dots&0 \end{pmatrix}^{t}$. Set $f_m(a):=\hat{a}$ and $f_m(s):=\hat{s}$. Set $g_0:=T_{1,-2}(s_0\hat{s}x)=T_{1,-2}(s_0\frac{s}{1}\frac{a}{s})=T_{1,-2}(s_0\hat{a})\in U_0$. By Lemma 4.5, $[^{\epsilon_0}h,g_0]\not\in U_{2n}((R_m,\Lambda_m),(I_m,$ $I_m\cap\Lambda_m))$. Since $U_{2n}((R_m,\Lambda_m),(I_m,$ $I_m\cap\Lambda_m))$ is normal, it follows that $\sigma:=$ $^{\epsilon_1}[^{\epsilon_0}h,g_0]\not\in U_{2n}((R_m,\Lambda_m),(I_m,I_m\cap\Lambda_m))$. Since 
\begin{align*}
&\sigma\\
=&^{\epsilon_1}[^{\epsilon_0}h,g_0]\\
=&^{\epsilon_1}(g_0^{-1}+(^{\epsilon_0} h)_{*1}s_0\hat{a}((^{\epsilon_0} h)^{-1})_{-2,*}g_0^{-1}\\
&-(^{\epsilon_0} h)_{*2}\bar\lambda_m\overline{s_0\hat{a}}((^{\epsilon_0} h)^{-1})_{-1,*}g_0^{-1})\\
=&^{\epsilon_1}(g_0^{-1})+^{\epsilon_1}((^{\epsilon_0} h)_{*1}s_0\hat{a}((^{\epsilon_0} h)^{-1})_{-2,*}g_0^{-1})\\
&-^{\epsilon_1}((^{\epsilon_0} h)_{*2}\bar\lambda_m\overline{s_0\hat{a}}((^{\epsilon_0} h)^{-1})_{-1,*}g_0^{-1})\\
=&^{\epsilon_1}(g_0^{-1})+\epsilon_1(^{\epsilon_0} h)_{*1}s_0\hat{a}((^{\epsilon_0} h)^{-1})_{-2,*}g_0^{-1}(\epsilon_1)^{-1}\\
&-\epsilon_1(^{\epsilon_0} h)_{*2}\bar\lambda_m\overline{s_0\hat{a}}((^{\epsilon_0} h)^{-1})_{-1,*})g_0^{-1}(\epsilon_1)^{-1}
\end{align*}
and $^{\epsilon_1}(g_0^{-1})=e+(\epsilon_1)_{*1}s_0\hat{a}((\epsilon_1)^{-1})_{-2,*}-(\epsilon_1)_{*2}\bar\lambda_m\overline{s_0\hat{a}}((\epsilon_1)^{-1})_{-1,*}$ has the form 
\[\begin{pmatrix}e^{n\times n} & * \\ 0 & e^{n\times n} \end{pmatrix},\]
$\sigma$ has the form 
\[\left(\begin{array}{ccc|ccc}
\alpha_1 &&\alpha_2 &\beta_1 &&\beta_2
\\ 0 && e^{2\times 2} &\beta_3&&\beta_4
\\\hline &\gamma&&&\delta&
\end{array}\right)=\begin{pmatrix}\alpha&\beta \\ \gamma&\delta\end{pmatrix}\]
where $\alpha_1, \beta_2\in M_{n-2}(R_m)$, $\alpha_2, \beta_1\in M_{(n-2)\times 2}(R_m)$, $\beta_3\in M_{2}(R_m)$, $\beta_4\in M_{2\times (n-2)}(R_m)$ and $\alpha=(a_{ij})_{1\leq i,j\leq n}, \beta=
(b_{ij})_{\substack{1 \leq i \leq n\\-n \leq j \leq -1}},\gamma=(c_{ij})_{\substack{-n\leq i \leq -1\\1 \leq j \leq n}},\delta=(d_{ij})_{-n \leq i,j \leq -1}\in M_n($ $R_m)$.\\
\\
\underline{case 1} Assume that either $\alpha\not\equiv e^{n\times n} (\textnormal{mod }I_m)$ or $\gamma\not\equiv 0 (\textnormal{mod }I_m)$ or $\delta\not\equiv e^{n\times n} (\textnormal{mod }$ $I_m)$.\\ We will show that it follows that $\alpha\not\equiv e^{n\times n} (\textnormal{mod }I_m)$ or $\gamma\not\equiv 0 (\textnormal{mod }I_m)$. Assume that $\alpha\equiv e^{n\times n} (\textnormal{mod }I_m)$ and $\gamma\equiv 0 (\textnormal{mod }I_m)$. Let $\kappa: M_{n}(R_m)\rightarrow M_{n}(R_m/I_m)$ be the homomorphism induced by the canonical homomorphism $R_m\rightarrow R_m/I_m$. Since the image of $\sigma$ in $U_{2n}(R_m/I_m,\Lambda_m/(\Lambda_m\cap I_m))$ equals $\begin{pmatrix}
e^{n\times n}& \kappa(\beta)\\ 0 &\kappa(\delta) \end{pmatrix}$, $\kappa(\delta)=e^{n\times n}$ by Lemma 3.9. That is equivalent to $\delta\equiv e^{n\times n}(\textnormal{mod }I_m)$. Since this is a contradiction, $\alpha\not\equiv e^{n\times n} (\textnormal{mod }I_m)$ or $\gamma\not\equiv 0 (\textnormal{mod }I_m)$. Hence there is an $i\in\{1,\dots,n\}$ such that $\sigma_{*i}\not\equiv e_i (\textnormal{mod }I_m)$.\\
\\
\underline{case 1.1} Assume that $i\in\{1,\dots,n-2\}$.\\Clearly the $(n-1)$-th row of
\small
\begin{align*}
&[\sigma,T_{i,-(n-1)}(1)]\\
=&(e+\sigma_{*i}\sigma'_{-(n-1),*}-\bar\lambda_m\sigma_{*,n-1}\sigma'_{-i,*})T_{i,-(n-1)}(-1)\\
=&(e+\bordermatrix{ & \cr 1&a_{1i} \cr &\vdots \cr &a_{n-2,i} \cr &0 \cr n&0 \cr -n&c_{-n,i}\cr &\vdots\cr -1&c_{-1,i}}
{\let\quad\thinspace \bordermatrix{&1&&n&-n&&&&-1\cr&\lambda_m\bar c_{-1,n-1}&\dots&\lambda_m\bar c_{-n,n-1} & 0 & 1~ & \bar a_{n-2,n-1}&\dots& \bar a_{1,n-1} }}
\end{align*}
\begin{align*}
&-\bar\lambda_m\bordermatrix{ & \cr 1&a_{1,n-1} \cr &\vdots \cr &a_{n-2,n-1}\cr &1 \cr n&0 \cr -n&c_{-n,n-1}\cr &\vdots\cr -1&c_{-1,n-1}}
{\let\quad\thinspace\bordermatrix{&1&&n&-n&&&&-1\cr&\lambda_m\bar c_{-1,i}&\dots&\lambda_m\bar c_{-n,i} & 0 & 0~ &\bar a_{n-2,i}&\dots& \bar a_{1i}})}\\
&\cdot T_{i,-(n-1)}(-1)
\end{align*}
\normalsize
is not congruent to $f_{n-1}$ modulo $I_m$ since $\sigma_{*i}\not\equiv e_i (\textnormal{mod }I_m)$. Recall that $f_l=e_l^t$ for any $l\in\{1,\dots,-1\}$. Hence $[\sigma,T_{i,-(n-1)}(1)]\not\in U_{2n}((R_m,\Lambda_m),(I_m,I_m\cap\Lambda_m))$. Set $g_1:=T_{i,-(n-1)}(s_1)\in U_1$. By Lemma 4.5, $[\sigma,g_1]\not\in U_{2n}((R_m,\Lambda_m),(I_m,I_m\cap\Lambda_m))$. Clearly the $n$-th row of $[\sigma,g_1]$ equals $f_n$. Let $P_{1n}$ be as in 3.11. Set $\epsilon_2:=P_{1n}\in EU_{2n}(R_m,\Lambda_m)$. Then the first row of $\tau:=$$^{\epsilon_2}[\sigma,g_1]$ equals $f_1$. Since $U_{2n}((R_m,\Lambda_m),(I_m,I_m\cap\Lambda_m))$ is normal,  $\tau \not\in U_{2n}((R_m,\Lambda_m),(I_m,I_
m\cap\Lambda_m))$. Clearly $\tau$ has the form 
 \[\left(\begin{array}{cc|cc}
 1 &0&0&0
 \\ A_3&A_4&B_3&0
 \\ \hline C_1 &C_2& D_1&0
 \\ C_3&C_4&D_3&1
 \end{array}\right)=\begin{pmatrix}A&B \\ C&D\end{pmatrix}\]
 where $C_3\in R_m$, $A_3, C_1\in (R_m)^{n-1}$, $C_4, D_3\in$ $^{n-1}(R_m)$, $A_4,B_3,C_2,D_1\in M_{n-1}(R_m)$ and $A=(A_{ij})_{1\leq i,j\leq n}, B=
 (B_{ij})_{\substack{1 \leq i \leq n\\-n \leq j \leq -1}},C=(C_{ij})_{\substack{-n\leq i \leq -1\\1 \leq j \leq n}},D=(D_{ij})_{-n \leq i,j \leq -1}\in M_n(R_m)$. Set 
 \[E=(E_{ij})_{2\leq i,j\leq -2}:=\left(\begin{array}{c|c}
A_4&B_3\\ \hline C_2&D_1
  \end{array}\right)\in M_{2n-2}(R_m).\]
 \underline{case 1.1.1} Assume that $E\not\equiv e^{(2n-2)\times(2n-2)} (\textnormal{mod }I_m)$.\\ There are $i,j\in\{2,\dots,$ $-2\}$ such that $(E-e^{(2n-2)\times(2n-2)})_{ij}\not\in I_m$. Set $g_2:=T_{-1,i}(s_2)\in U_2$. 
 Then $\omega:=[\tau^{-1},g_2]$ has the form 
 \[\begin{pmatrix} 1 &0&0 \\ *&e^{(2n-2)\times (2n-2)}& 0\\ *& w&1\end{pmatrix}\]
 where $w=(w_2,\dots,w_{-2})=s_2(E-e^{(2n-2)\times (2n-2)})_{i*}$. Since $(E-e^{(2n-2)\times(2n-2)})_{ij}$ $\not\in I_m$, $w_j=:\frac{b'}{t'}\not\in I_m$. Set $b:=\frac{b'}{1}\in R_m$ and  $t:=\frac{t'}{1}\in R_m$. Choose an $l\neq\pm 1, \pm j$ and set $g_3:=T_{jl}(s_3)\in U_3$, $g_4:=T_{l,-j}(s_4s_5t)\in U_4$ and $g_5:=T_{-1,-j}(s_3s_4s_5tw_j)=T_{-1,-j}(s_3s_4s_5b)\in U_5$. Notice that $g_5\not \in U_{2n}((R_m,\Lambda_m),(I_m,I_m\cap\Lambda_m))$, since $w_j\not\in I_m$ and $s_3s_4s_5t$ is invertible. One checks easily that 
 \begin{align*}
 &[[[^{\epsilon_2}([^{\epsilon_1}[^{\epsilon_0}h,g_0],g_2]^{-1}),g_2],g_3],g_4]\\
 =&[[[(^{\epsilon_2}[^{\epsilon_1}[^{\epsilon_0}h,g_0],g_2])^{-1},g_2],g_3],g_4]\\
 =&[[\omega,g_3],g_4]\\
 =&g_5.
 \end{align*}
 Set $g'_i:=\psi_m(g_i)~\forall i\in\{0,\dots,5\}$ and $\epsilon'_i:=\psi_m(\epsilon_i)~\forall i\in\{0,1,2\}$. Then 
 \[[[[^{\epsilon'_2}([^{\epsilon'_1}[^{\epsilon'_0}h',g'_0],g'_1]^{-1}),g'_2],g'_3],g'_4]=g'_5.\]
 \underline{case 1.1.2} Assume that $E\equiv e^{(2n-2)\times(2n-2)} (\textnormal{mod }I_m)$ and $A_3\equiv 0(\textnormal{mod }I_m)$.\\ Set $\xi_1:=\prod\limits_{l=2}^{n} T_{l1}(-A_{l1})\in EU_{2n}((R_m,\Lambda_m),(I_m,\Gamma_m))$. Since $\xi_1\in EU_{2n}((R_m,\Lambda_m),(I_m,$ $\Gamma_m))$ $\subseteq U_{2n}((R_m,\Lambda_m),(I_m,\Gamma_m))\subseteq U_{2n}((R_m,\Lambda_m),(I_m,I_m\cap\Lambda_m))$ and $\tau\not\in U_{2n}((R_m,$ $\Lambda_m),(I_m,I_m\cap\Lambda_m))$, $\xi_1\tau\not\in U_{2n}((R_m,\Lambda_m),(I_m,I_m\cap\Lambda_m))$. Clearly 
  \[\xi_1\tau=\left(\begin{array}{cc|cc}
  1 &0&0&0
  \\ 0&A_4&B_3&0
  \\ \hline C_1 &C_2& D_1&0
  \\ C'_3&C'_4&D'_3&1
  \end{array}\right)\]
  for some $C'_3\in R_m$ and $C'_4, D'_3\in$ $^{n-1}(R_m)$ such that $D'_3\equiv 0(\textnormal{mod }I_m)$ (consider the image of $\xi_1\tau$ in $U_{2n}(R_m/I_m,\Lambda_m/(\Lambda_m\cap I_m))$). \\
 \\
 \underline{case 1.1.2.1} Assume that there is an $i\in\{3,\dots,n\}$ such that $C'_{-1,i}\not\in I_m$.\\ Set $\epsilon_{21}:=T_{12}(-1)\in EU_{2n}(R_m,\Lambda_m)$. 
 Then $^{\epsilon_{21}}(\xi_1\tau)$ has the form 
   \[\left(\begin{array}{cc|cc}
   1 &A''_2&B''_1&B''_2
   \\ 0&A_4&B''_3&B''_4
   \\ \hline C''_1 &C''_2& D''_1&D''_2
   \\ C''_3&C''_4&D''_3&D''_4
   \end{array}\right)\]
   where $B''_2,C''_3,D''_4\in R_m$, $B''_4,C''_1\in (R_m)^{n-1}$, $A''_2,B''_1,C''_4, D''_3\in$$^{n-1}(R_m)$, $B''_3,C''_2,D''_1\in M_{n-1}(R_m)$. Furthermore $A''_2\equiv 0(\textnormal{mod }I_m)$ and $C''_{-2,i}\equiv C'_{-1,i}(\textnormal{mod }I_m)$. Set $\xi_2:=\prod\limits_{l=2}^{n} T_{1l}(-A''_{1l})\in EU_{2n}((R_m,\Lambda_m),(I_m,\Gamma_m))$. $\omega:=$$^{\epsilon_{21}}(\xi_1\tau)\xi_2$ has the form 
   \[\left(\begin{array}{cc|cc}
      1 &0&B'''_1&B'''_2
      \\ 0&A_4&B'''_3&B'''_4
      \\ \hline C'''_1 &C'''_2& D'''_1&D'''_2
      \\ C'''_3&C'''_4&D'''_3&D'''_4
      \end{array}\right)\]
      where $B'''_2,C'''_3, D'''_4\in R_m$, $B'''_4,C'''_1, D'''_2\in (R_m)^{n-1}$, $B'''_1,C'''_4, D'''_3\in$$^{n-1}(R_m)$, $B'''_3,C'''_2,$ $D'''_1\in M_{n-1}(R_m)$.  Further $C'''_{-2,i}\equiv C''_{-2,i}\equiv C'_{-1,i}(\textnormal{mod }I_m)$. Since $C'_{-1,i}\not\in I_m$, $C'''_{-2,i}\not\in I_m$. Set $g_2:=T_{2,-i}(s_2)\in U_2$. Then
     \begin{align*}
      &[\omega,g_2]\\
      =&(e+\omega_{*2}s_2\omega'_{-i,*}-\omega_{*i}\bar\lambda_m s_2\omega'_{-2,*})(g_2)^{-1}\\
      =&(e+s_2\bordermatrix{ & \cr 1& 0 \cr &A_{22} \cr & \vdots \cr n&A_{n2} \cr -n&C'''_{-n,2} \cr &\vdots \cr -1&C'''_{-1,2}}
      \bordermatrix{&1&&n&-n&&&-1\cr&\lambda_m\bar C'''_{-1,i}&\dots&\lambda_m\bar C'''_{-n,i} & \bar A_{ni} & \dots& \bar A_{2i}& 0}\\
      &-\bar\lambda_m s_2\bordermatrix{ & \cr 1& 0 \cr &A_{2i} \cr & \vdots \cr n&A_{ni} \cr -n&C'''_{-n,i} \cr &\vdots \cr -1&C'''_{-1,i}}
            \bordermatrix{&1&&n&-n&&&-1\cr&\lambda_m\bar C'''_{-1,2}&\dots&\lambda_m\bar C'''_{-n,2} & \bar A_{n2} & \dots& \bar A_{22}& 0})\\
      &\cdot (g_2)^{-1}.
      \end{align*}
      Clearly $[\omega,g_2]_{1*}=f_1$ and 
      \[
      [\omega,g_2]_{22}=1+s_2A_{22}\lambda_m \bar C'''_{-2,i}-\bar\lambda_ms_2A_{2i}\lambda_m \bar C'''_{-2,2}.
      \]
      Since $A_{2i}\in I_m$, $-\bar\lambda_ms_2A_{2i}\lambda_m \bar C'''_{-2,2}\in I_m$. Since $C'''_{-2,i}\not\in I_m$, $\bar C'''_{-2,i}\not\in I_m$. Hence $s_2A_{22}\lambda_m \bar{C'''}_{-2,i}\not\in I_m$ since $A_{22}\equiv 1(\textnormal{mod }I_m)$ and $s_2$ and $\lambda_m$ are invertible. It follows that $[\omega,g_2]_{22}\not\equiv 1(\textnormal{mod }I_m)$. One can proceed now as in case 1.1.1 (note that $\psi_m(\xi_1)=\psi_m(\xi_2)=e$ since $\xi_1,\xi_2\in EU_{2n}((R_m,\Lambda_m),(I_m,\Gamma_m))\subseteq U_{2n}((R_m,\Lambda_m),(I_m,$ $\Gamma_m)))$.\\ 
\\
\underline{case 1.1.2.2} Assume that $C'_{-1,2}\not\in I_m$.\\ 
Set $\epsilon_{21}:=T_{13}(-1)\in EU_{2n}(R_m,\Lambda_m)$.
Then $^{\epsilon_{21}}(\xi_1\tau)$ has the form 
   \[\left(\begin{array}{cc|cc}
   1 &A''_2&B''_1&B''_2
   \\ 0&A_4&B''_3&B''_4
   \\ \hline C''_1 &C''_2& D''_1&D''_2
   \\ C''_3&C''_4&D''_3&D''_4
   \end{array}\right)\]
   where $B''_2,C''_3,D''_4\in R_m$, $B''_4,C''_1\in (R_m)^{n-1}$, $A''_2,B''_1,C''_4, D''_3\in$$^{n-1}(R_m)$, $B''_3,C''_2,D''_1\in M_{n-1}(R_m)$. Furthermore $A''_2\equiv 0(\textnormal{mod }I_m)$ and $C''_{-3,2}\equiv C'_{-1,2}(\textnormal{mod }I_m)$. Set $\xi_2:=\prod\limits_{l=2}^{n} T_{1l}(-A''_{1l})\in EU_{2n}((R_m,\Lambda_m),(I_m,\Gamma_m))$. Then $\omega:=$$^{\epsilon_{21}}(\xi_1\tau)\xi_2$ has the form 
   \[\left(\begin{array}{cc|cc}
      1 &0&B'''_1&B'''_2
      \\ 0&A_4&B'''_3&B'''_4
      \\ \hline C'''_1 &C'''_2& D'''_1&D'''_2
      \\ C'''_3&C'''_4&D'''_3&D'''_4
      \end{array}\right)\]
      where $B'''_2,C'''_3, D'''_4\in R_m$, $B'''_4,C'''_1, D'''_2\in (R_m)^{n-1}$, $B'''_1,C'''_4, D'''_3\in$$^{n-1}(R_m)$, $B'''_3,C'''_2,$ $D'''_1\in M_{n-1}(R_m)$.  Further $C'''_{-3,2}\equiv C''_{-3,2}\equiv C'_{-1,2}(\textnormal{mod }I_m)$. Since $C'_{-1,2}\not\in I_m$, $C'''_{-3,2}\not\in I_m$. Set $g_2:=T_{3,-2}(s_2)\in U_2$. Then
      \begin{align*}
      &[\omega,g_2]\\
      =&(e+\omega_{*3}s_2\omega'_{-2,*}-\omega_{*2}\bar\lambda_ms_2\omega'_{-3,*})(g_2)^{-1}\\
      =&(e+s_2\bordermatrix{ & \cr 1& 0 \cr &A_{23} \cr & \vdots \cr n&A_{n3} \cr -n&C'''_{-n,3} \cr &\vdots \cr -1&C'''_{-1,3}}
      \bordermatrix{&1&&n&-n&&&-1\cr&\lambda_m\bar C'''_{-1,2}&\dots&\lambda_m\bar C'''_{-n,2} & \bar A_{n2} & \dots& \bar A_{22}& 0}\\
      &-\bar\lambda_ms_2\bordermatrix{ & \cr 1& 0 \cr &A_{22} \cr & \vdots \cr n&A_{n2} \cr -n&C'''_{-n,2} \cr &\vdots \cr -1&C'''_{-1,2}}
      \bordermatrix{&1&&n&-n&&&-1\cr&\lambda_m\bar C'''_{-1,3}&\dots&\lambda_m\bar C'''_{-n,3} & \bar A_{n3} & \dots& \bar A_{23}& 0})
            \\
      &\cdot (g_2)^{-1}.
      \end{align*}
      Clearly $[\omega,g_2]_{1*}=f_1$ and 
      \[
      [\omega,g_2]_{33}=1+s_2A_{33}\lambda_m \bar C'''_{-3,2}-\bar\lambda_ms_2A_{32}\lambda_m \bar C'''_{-3,3}.
      \]
      Since $A_{32}\in I_m$, $-\bar\lambda_ms_2A_{32}\lambda_m \bar C'''_{-3,3}\in I_m$. Since $C'''_{-3,2}\not\in I_m$, $\bar C'''_{-3,2}\not\in I_m$. Hence $s_2A_{33}\lambda_m \bar{C'''}_{-3,2}\not\in I_m$ since $A_{33}\equiv 1(\textnormal{mod }I_m)$ and $s_2$ and $\lambda_m$ are invertible. It follows that $[\omega,g_2]_{33}\not\equiv 1(\textnormal{mod }I_m)$. One can proceed now as in case 1.1.1\\ 
\\
\underline{case 1.1.2.3} Assume that $C'_{-1,i}\in I_m~\forall i\in\{2,\dots,n\}$.\\ It follows that $C_{i1}\in I_m~\forall i\in\{-n,\dots,-2\}$. Since $\xi_1\tau\not\in U_{2n}((R_m,\Lambda_m),(I_m,I_m\cap\Lambda_m))$, $C'_{-1,1}\not\in I_m$. Set $\epsilon_{21}:=T_{12}(-1)\in EU_{2n}(R_m,\Lambda_m)$. 
 Then $^{\epsilon_{21}}(\xi_1\tau)$ has the form 
   \[\left(\begin{array}{cc|cc}
   1 &A''_2&B''_1&B''_2
   \\ 0&A_4&B''_3&B''_4
   \\ \hline C''_1 &C''_2& D''_1&D''_2
   \\ C''_3&C''_4&D''_3&D''_4
   \end{array}\right)\]
   where $B''_2,C''_3,D''_4\in R_m$, $B''_4,C''_1\in (R_m)^{n-1}$, $A''_2,B''_1,C''_4, D''_3\in$$^{n-1}(R_m)$, $B''_3,C''_2,D''_1\in M_{n-1}(R_m)$. Furthermore $A''_2\equiv 0(\textnormal{mod }I_m)$ and $C''_{-2,2}=C_{-2,2}+C'_{-1,2}+C_{-2,1}+C'_{-1,1}\equiv C'_{-1,1}(\textnormal{mod }I_m)$. Set $\xi_2:=\prod\limits_{l=2}^{n} T_{1l}(-A''_{1l})\in EU_{2n}((R_m,\Lambda_m),(I_m,\Gamma_m))$. Then $\omega:=$$^{\epsilon_{21}}(\xi_1\tau)\xi_2$ has the form 
   \[\left(\begin{array}{cc|cc}
      1 &0&B'''_1&B'''_2
      \\ 0&A_4&B'''_3&B'''_4
      \\ \hline C'''_1 &C'''_2& D'''_1&D'''_2
      \\ C'''_3&C'''_4&D'''_3&D'''_4
      \end{array}\right)\]
      where $B'''_2,C'''_3, D'''_4\in R_m$, $B'''_4,C'''_1, D'''_2\in (R_m)^{n-1}$, $B'''_1,C'''_4, D'''_3\in$$^{n-1}(R_m)$, $B'''_3,C'''_2,$ $D'''_1\in M_{n-1}(R_m)$.  Further $C'''_{-2,2}\equiv C''_{-2,2}\equiv C'_{-1,1}(\textnormal{mod }I_m)$. Since $C'_{-1,1}\not\in I_m$, $C'''_{-2,2}\not\in I_m$. Set $g_2:=T_{2,-3}(f_m(t_3s_0))\in U_2$. Then
      \begin{align*}
      &[\omega,g_2]\\
      =&(e+\omega_{*2}s_2\omega'_{-3,*}-\omega_{*3}\bar\lambda_ms_2\omega'_{-2,*})(g_2)^{-1}\\
      =&(e+s_2\bordermatrix{ & \cr 1& 0 \cr &A_{22} \cr & \vdots \cr n&A_{n2} \cr -n&C'''_{-n,2} \cr &\vdots \cr -1&C'''_{-1,2}}
      \bordermatrix{&1&&n&-n&&&-1\cr&\lambda_m\bar C'''_{-1,3}&\dots&\lambda_m\bar C'''_{-n,3} & \bar A_{n3} & \dots& \bar A_{23}& 0}\\
      &-\bar\lambda_ms_2\bordermatrix{ & \cr 1& 0 \cr &A_{23} \cr & \vdots \cr n&A_{n3} \cr -n&C'''_{-n,3} \cr &\vdots \cr -1&C'''_{-1,3}}
            \bordermatrix{&1&&n&-n&&&-1\cr&\lambda_m\bar C'''_{-1,2}&\dots&\lambda_m\bar C'''_{-n,2} & \bar A_{n2} & \dots& \bar A_{22}& 0})\\
      &\cdot (g_2)^{-1}.
      \end{align*}
      Clearly $[\omega,g_2]_{1*}=f_1$ and 
      \[
      [\omega,g_2]_{32}=s_2A_{32}\lambda_m \bar C'''_{-2,3}-\bar\lambda_ms_2A_{33}\lambda_m \bar C'''_{-2,2}.
      \]
      Since $A_{32}\in I_m$, $-\bar\lambda_ms_2A_{32}\lambda_m \bar C'''_{-2,2}\in I_m$. Since $C'''_{-2,2}\not\in I_m$, $\bar C'''_{-2,2}\not\in I_m$. Hence $s_2A_{33}\lambda_m \bar{C'''}_{-2,i}\not\in I_m$ since $A_{33}\equiv 1(\textnormal{mod }I_m)$ and $s_2$ and $\lambda_m$ are invertible. It follows that $[\omega,g_2]_{32}\not\in I_m$. One can proceed now as in case 1.1.1\\
\\
\underline{case 1.1.3} Assume that $E\equiv e^{(2n-2)\times(2n-2)} (\textnormal{mod }I_m)$ and $A_3\not\equiv 0(\textnormal{mod }I_m)$.\\ Since $A_3\not\equiv 0(\textnormal{mod }I_m)$, $D_3\not\equiv 0(\textnormal{mod }I_m)$. Hence there is an $i\in\{-n,\dots,-2\}$ such that $D_{-1,i}\not\in I_m$. Choose a $j\in \{2,\dots, n\}\setminus\{-i\}$ and set $g_2:=T_{ij}(s_2)\in U_2$. Then
\small
\begin{align*}
&[\tau,g_2]\\
      =&(e+\tau_{*i}s_2\tau'_{j*}-\tau_{*,-j}\lambda_ms_2\tau'_{-i,*})(g_2)^{-1}\\
      =&(e+s_2\bordermatrix{ & \cr 1& 0 \cr &B_{2i} \cr & \vdots \cr n&B_{ni} \cr -n&D_{-n,i} \cr &\vdots \cr -1&D_{-1,i}}
      \bordermatrix{&1&&n&-n&&&-1\cr&\bar D_{-1,-j}&\dots&\bar D_{-n,-j} & \bar\lambda_m\bar B_{n,-j} & \dots& \bar\lambda_m\bar B_{2,-j}& 0}\\
      &-\lambda_ms_2\bordermatrix{ & \cr 1& 0 \cr &B_{2,-j} \cr & \vdots \cr n&B_{n,-j} \cr -n&D_{-n,-j} \cr &\vdots \cr -1&D_{-1,-j}}
            \bordermatrix{&1&&n&-n&&&-1\cr&\bar D_{-1,i}&\dots&\bar D_{-n,i} & \bar\lambda_m\bar B_{ni} & \dots& \bar\lambda_m\bar B_{2i}& 0})\\
      &\cdot (g_2)^{-1}.\hspace{11cm}\hbox{}
      \end{align*}
      \normalsize
      Clearly $[\tau,g_2]_{1*}=f_1$ and 
      \begin{align*}
      &[\tau,g_2]_{-1,j}\\
      =&s_2D_{-1,i}\bar D_{-j,-j}-\lambda_ms_2D_{-1,-j}\bar D_{-j,i}\\
      &-s_2(s_2D_{-1,i}\bar \lambda_m\bar B_{-i,-j}-\lambda_ms_2D_{-1,-j}\bar \lambda_m\bar B_{-i,i})
      \end{align*}
      Since $\bar D_{-j,i},\bar B_{-i,-j}, \bar B_{-i,i}\in I_m$, it follows that $-\lambda_ms_2D_{-1,-j}\bar D_{-j,i}-s_2(s_2D_{-1,i}\bar \lambda_m$ $\bar B_{-i,-j}-\lambda_ms_2$ $D_{-1,-j}\bar \lambda_m\bar B_{-i,i})\in I_m$. On the other hand $s_2D_{-1,i}\bar D_{-j,-j}\not\in I_m$ since $D_{1,-i}\not\in I_m$, $\bar D_{-j,-j}\equiv 1(\textnormal{mod }I_m)$ and $s_2$ is invertible. It follows that $[\tau,g_2]_{-1,j}\not\in I_m$ and hence $[\tau,g_2]\not\in U_{2n}((R_m,\Lambda_m),(I_m,I_m\cap\Lambda_m))$. Further $[\tau,g_2]_{l1}\in I_m~\forall l\in\{2,\dots,n\}$ since $B\equiv 0(\textnormal{mod }I_m)$. Thus one can proceed now as in case 1.1.1 or as in case 1.1.2.\\
\\
\underline{case 1.2} Assume $\sigma_{*j}\equiv e_j(\textnormal{mod }I_m)~\forall j\in\{1,\dots,n-2\}$, $\sigma_{*,n-1}\not\equiv e^{n-1}(\textnormal{mod }I_m)$ and $a_{1,n-1}\in I_m$.\\
Consider the first row of 
\small
\begin{align*}
&[\sigma,T_{1,-(n-1)}(1)]\\
=&(e+\sigma_{*1}\sigma'_{-(n-1),*}-\bar\lambda_m\sigma_{*,n-1}\sigma'_{-1,*})T_{1,-(n-1)}(-1)\\
=&(e+\bordermatrix{ & \cr 1&a_{11} \cr &\vdots \cr &a_{n-2,1}\cr & 0\cr n&0 \cr -n&c_{-n,1}\cr &\vdots\cr -1&c_{-1,1}}
{\let\quad\thinspace \bordermatrix{&1&&n&-n&&&&-1\cr&\lambda_m\bar c_{-1,n-1}&\dots&\lambda_m\bar c_{-n,n-1} & 0 & 1 ~& \bar a_{n-2,n-1} & \dots& \bar a_{1,n-1} }}\\
&-\bar\lambda_m\bordermatrix{ & \cr 1&a_{1,n-1} \cr &\vdots \cr &a_{n-2,n-1} \cr &1 \cr n&0 \cr -n&c_{-n,n-1}\cr &\vdots\cr -1&c_{-1,n-1}}
{\let\quad\thinspace\bordermatrix{&1&&n&-n&&&&-1\cr&\lambda_m\bar c_{-1,1}&\dots&\lambda_m\bar c_{-n,1} & 0  &0~ & \bar a_{n-2,1}&\dots& \bar a_{11}})}\\
&\cdot T_{1,-(n-1)}(-1)
\end{align*}
\normalsize
which equals 
\setlength\jot{0.5cm}
\begin{align*}
&\bordermatrix{&1&&&n&-n&&-1\cr&1&0 &\dots& 0 & 0& \dots&0 }\\
&+a_{11}\bordermatrix{&1&&n&-n&&&&-1\cr&\lambda_m\bar c_{-1,n-1}&\dots&\lambda_m\bar c_{-n,n-1} & 0 & 1 ~& \bar a_{n-2,n-1} & \dots& \bar a_{1,n-1,} }
\end{align*}
\begin{align*}
&-\bar\lambda_m a_{1,n-1}\bordermatrix{&1&&n&-n&&&&-1\cr&\lambda_m\bar c_{-1,1}&\dots&\lambda_m\bar c_{-n,1} & 0  &0~ & \bar a_{n-2,1}&\dots& \bar a_{11}}\\
&+\bordermatrix{&1&&n&-n&&&&&-1\cr&0&\dots&0&0&x_1&0&\dots&0&x_2}
\end{align*}\setlength\jot{3pt}for some $x_1,x_2\in R$. It is clearly not
congruent to $f_1$ modulo $I_m$ since $a_{11}\equiv 1(\textnormal{mod }I_m)$, $\sigma_{*,n-1}\not\equiv e^{n-1} (\textnormal{mod }I_m)$ and $a_{1,n-1}\in I_m$. Hence $[\sigma,T_{1,-(n-1)}(1)]\not\in U_{2n}((R_m,\Lambda_m),(I_m,I_m\cap\Lambda_m))$. Set $g_1:=T_{1,-(n-1)}(s_1)\in U_1$. By Lemma 4.5, $[\sigma,g_1]\not\in U_{2n}((R_m,\Lambda_m),(I_m,I_m\cap\Lambda_m))$. Clearly the $n$-th row of $[\sigma,g_1]$ equals $f_n$. Let $P_{1n}$ be as in 3.11 and set $\epsilon_2:=P_{1n}\in EU_{2n}(R_m,\Lambda_m)$. Then the first row of $\tau:=$$^{\epsilon_2}[\sigma,g_1]$ equals $f_1$. Since $U_{2n}((R_m,\Lambda_m),(I_m,I_m\cap\Lambda_m))$ is normal,  $\tau \not\in U_{2n}((R_m,\Lambda_m),(I_m,I_m\cap\Lambda_m))$. One can proceed now as in case 1.1.\\
\\
\underline{case 1.3} Assume $\sigma_{*j}\equiv e_j(\textnormal{mod }I_m)~\forall j\in\{1,\dots,n-2\}$ and $a_{1,n-1}\not\in I_m$.\\
Consider the second row of 
\small
\begin{align*}
&[\sigma,T_{2,-(n-1)}(1)]\\
=&(e+\sigma_{*2}\sigma'_{-(n-1),*}-\bar\lambda_m\sigma_{*,n-1}\sigma'_{-2,*})T_{1,-(n-1)}(-1)\\
=&(e+\bordermatrix{ & \cr 1&a_{12} \cr &\vdots \cr &a_{n-2,2} \cr &0 \cr n&0 \cr -n&c_{-n,2}\cr &\vdots\cr -1&c_{-1,2}}
{\let\quad\thinspace \bordermatrix{&1&&n&-n&&&&-1\cr&\lambda_m\bar c_{-1,n-1}&\dots&\lambda_m\bar c_{-n,n-1} & 0 & 1~ & \bar a_{n-2,n-1}& \dots&\bar a_{1,n-1} }}\\
&-\bar\lambda_m\bordermatrix{ & \cr 1&a_{1,n-1} \cr &\vdots\cr &a_{n-2,n-1} \cr &1 \cr n&0 \cr -n&c_{-n,n-1}\cr &\vdots\cr -1&c_{-1,n-1}}
{\let\quad\thinspace\bordermatrix{&1&&n&-n&&&&-1\cr&\lambda_m\bar c_{-1,2}&\dots&\lambda_m\bar c_{-n,2} & 0 &0~& \bar a_{n-2,2}&\dots& \bar a_{12}})}\\
&\cdot T_{2,-(n-1)}(-1)
\end{align*}
\normalsize
which equals 
\setlength\jot{0.4cm}
\begin{align*}
&\bordermatrix{&1&&&&n&-n&&-1\cr&0&1&0 &\dots& 0 & 0& \dots&0 }\\
&+a_{22}\bordermatrix{&1&&n&-n&&&&-1\cr&\lambda_m\bar c_{-1,n-1}&\dots&\lambda_m\bar c_{-n,n-1} & 0 & 1~ & \bar a_{n-2,n-1}& \dots&\bar a_{1,n-1} }\\
&-\bar\lambda_m a_{2,n-1}\bordermatrix{&1&&n&-n&&&&-1\cr&\lambda_m\bar c_{-1,2}&\dots&\lambda_m\bar c_{-n,2} & 0 &0~& \bar a_{n-2,2}&\dots& \bar a_{12}}\\
&+\bordermatrix{&1&&n&-n&&&&&&-1\cr&0&\dots&0&0&x_1&0&\dots&0&x_2&0}
\end{align*}
\setlength\jot{3pt}for some $x_1,x_2\in R$. Its last entry clearly does not lie in $I_m$. Hence $[\sigma,T_{2,-(n-1)}(1)]$ $\not\in U_{2n}((R_m,\Lambda_m),(I_m,I_m\cap\Lambda_m))$. Set $g_1:=T_{2,-(n-1)}(s_1)\in U_1$. By Lemma 4.5, $[\sigma,g_1]\not\in U_{2n}((R_m,\Lambda_m),(I_m,I_m\cap\Lambda_m))$. Clearly the $n$-th row of $[\sigma,g_1]$ equals $f_n$. Set $\epsilon_2:=P_{1n}\in EU_{2n}(R_m,\Lambda_m)$. Then the first row of $\tau:=$$^{\epsilon_2}[\sigma,g_1]$ equals $f_1$. Since $U_{2n}((R_m,\Lambda_m),(I_m,I_m\cap\Lambda_m))$ is normal,  $\tau \not\in U_{2n}((R_m,\Lambda_m),(I_m,I_m\cap\Lambda_m))$. One can proceed now as in case 1.1.\\
\\
\underline{case 1.4} Assume $\sigma_{*j}\equiv e_j(\textnormal{mod }I_m)~\forall j\in\{1,\dots,n-1\}$, $\sigma_{*n}\not\equiv e_n(\textnormal{mod }I_m)$ and $a_{1n}\in I_m$.\\
Consider the first row of 
\small
\begin{align*}
&[\sigma,T_{1,-n}(1)]\\
=&(e+\sigma_{*1}\sigma'_{-n,*}-\bar\lambda_m\sigma_{*n}\sigma'_{-1,*})T_{1,-n}(-1)\\
=&(e+\bordermatrix{ & \cr 1&a_{11}  \cr &\vdots \cr&a_{n-2,1}\cr &0\cr n&0 \cr -n&c_{-n,1}\cr &\vdots\cr -1&c_{-1,1}}
{\let\quad\thinspace \bordermatrix{&1&&n&-n&&&&-1\cr&\lambda_m\bar c_{-1,n}&\dots&\lambda_m\bar c_{-n,n} &1& 0~ &\bar a_{n-2,n}&\dots& \bar a_{1n} }}
\end{align*}
\begin{align*}
&-\bar\lambda_m\bordermatrix{ & \cr 1&a_{1n} \cr&\vdots \cr &a_{n-2,n}&\cr &0  \cr n&1 \cr -n&c_{-n,n}\cr &\vdots\cr -1&c_{-1,n}}
{\let\quad\thinspace\bordermatrix{&1&&n&-n&&&&-1\cr&\lambda_m\bar c_{-1,1}&\dots&\lambda_m\bar c_{-n,1} & 0  &0~& \bar a_{n-2,1}&\dots& \bar a_{11}})}\\
&\cdot T_{1,-n}(-1)
\end{align*}
\normalsize
which equals \begin{align*}
&\bordermatrix{&1&&&n&-n&&-1\cr&1&0 &\dots& 0 & 0& \dots&0 }\\
&+a_{11}\bordermatrix{&1&&n&-n&&&&-1\cr&\lambda_m\bar c_{-1,n}&\dots&\lambda_m\bar c_{-n,n} &1& 0~ &\bar a_{n-2,n}&\dots& \bar a_{1n} }\\
&-\bar\lambda_m a_{1n}\bordermatrix{&1&&n&-n&&&&-1\cr&\lambda_m\bar c_{-1,1}&\dots&\lambda_m\bar c_{-n,1} & 0  &0~& \bar a_{n-2,1}&\dots& \bar a_{11}}\\
&+\bordermatrix{&1&&n&-n&&&&-1\cr&0&\dots&0&x_1&0&\dots&0&x_2}
\end{align*}
for some $x_1,x_2\in R$. It is clearly not congruent to $f_1$ modulo $I_m$ since $a_{11}\equiv 1(\textnormal{mod }I_m)$, $\sigma_{*n}\not\equiv e_n (\textnormal{mod }I_m)$ and $a_{1n}\in I_m$. Hence $[\sigma,T_{1,-n}(1)]\not\in U_{2n}((R_m,\Lambda_m),$ $(I_m,I_m\cap\Lambda_m))$. Set $g_1:=T_{1,-n}(s_1)\in U_1$. By Lemma 4.5, $[\sigma,g_1]\not\in U_{2n}((R_m,\Lambda_m),$ $(I_m,I_m\cap\Lambda_m))$. Clearly the $(n-1)$-th row of $[\sigma,g_1]$ equals $f_{n-1}$. Set $\epsilon_2:=P_{1(n-1)}\in EU_{2n}(R_m,\Lambda_m)$. Then the first row of $\tau:=$$^{\epsilon_2}[\sigma,g_1]$ equals $f_1$. Since $U_{2n}((R_m,\Lambda_m),$ $(I_m,I_m\cap\Lambda_m))$ is normal,  $\tau \not\in U_{2n}((R_m,\Lambda_m),(I_m,I_m\cap\Lambda_m))$. One can proceed now as in case 1.1.\\
\\
\underline{case 1.5} Assume $\sigma_{*j}\equiv e_j(\textnormal{mod }I_m)~\forall j\in\{1,\dots,n-1\}$ and $a_{1n}\not\in I_m$.\\
Consider the second row of 
\small
\begin{align*}
&[\sigma,T_{2,-n}(1)]\\
=&(e+\sigma_{*2}\sigma'_{-n,*}-\bar\lambda_m\sigma_{*n}\sigma'_{-2,*})T_{2,-n}(-1)
\end{align*}
\begin{align*}
=&(e+\bordermatrix{ & \cr 1&a_{12} \cr&\vdots\cr &a_{n-2,2} \cr &0 \cr n&0 \cr -n&c_{-n,2}\cr &\vdots\cr -1&c_{-1,2}}
{\let\quad\thinspace \bordermatrix{&1&&n&-n&&&&-1\cr&\lambda_m\bar c_{-1,n}&\dots&\lambda_m\bar c_{-n,n} &1& 0~&\bar a_{n-2,n}&\dots & \bar a_{1n} }}\\
&-\bar\lambda_m\bordermatrix{ & \cr 1&a_{1n} \cr &\vdots\cr&a_{n-2,n} \cr &0  \cr n&1 \cr -n&c_{-n,n}\cr &\vdots\cr -1&c_{-1,n}}
{\let\quad\thinspace\bordermatrix{&1&&n&-n&&&&-1\cr&\lambda_m\bar c_{-1,2}&\dots&\lambda_m\bar c_{-n,2} & 0  &0~& \bar a_{n-2,2}&\dots& \bar a_{12}})}\\
&\cdot T_{2,-n}(-1)
\end{align*}
\normalsize
which equals \begin{align*}
&\bordermatrix{&1&&&&n&-n&&-1\cr&0&1&0 &\dots& 0 & 0& \dots&0 }\\
&+a_{22}\bordermatrix{&1&&n&-n&&&&-1\cr&\lambda_m\bar c_{-1,n}&\dots&\lambda_m\bar c_{-n,n} &1& 0~&\bar a_{n-2,n}&\dots & \bar a_{1n} }\\
&-\bar\lambda_m a_{2n}\bordermatrix{&1&&n&-n&&&&-1\cr&\lambda_m\bar c_{-1,2}&\dots&\lambda_m\bar c_{-n,2} & 0  &0~& \bar a_{n-2,2}&\dots& \bar a_{12}}\\
&+\bordermatrix{&1&&n&-n&&&&&-1\cr&0&\dots&0&x_1&0&\dots&0&x_2&0}
\end{align*}
for some $x_1,x_2\in R$. Its last entry does clearly not lie in
$I_m$ and hence $[\sigma,T_{2,-n}(1)]\not\in U_{2n}((R_m,$ $\Lambda_m),(I_m,I_m\cap\Lambda_m))$. Set $g_1:=T_{2,-n}(s_1)\in U_1$. By Lemma 4.5, $[\sigma,g_1]\not\in U_{2n}((R_m,$ $\Lambda_m),(I_m,I_m\cap\Lambda_m))$. Clearly the $(n-1)$-th row of $[\sigma,g_1]$ equals $f_{n-1}$. Set $\epsilon_2:=P_{1(n-1)}\in EU_{2n}(R_m,\Lambda_m)$. Then the first row of $\tau:=$$^{\epsilon_2}[\sigma,g_1]$ equals $f_1$. Since $U_{2n}((R_m,\Lambda_m),$ $(I_m,I_m\cap\Lambda_m))$ is normal,  $\tau \not\in U_{2n}((R_m,\Lambda_m),(I_m,I_m\cap\Lambda_m))$. One can proceed now as in case 1.1.\\
\\
\underline{case 2} Assume that $\alpha, \delta\equiv e^{n\times n} (\textnormal{mod }I_m)$ and $\gamma\equiv 0 (\textnormal{mod }I_m)$.\\ Recall that $\sigma=$$^{\epsilon_1}[^{\epsilon_0}h,$ $g_0]\not\in U_{2n}((R_m,\Lambda_m),(I_m,I_m\cap\Lambda_m))$ has the form 
      \[\left(\begin{array}{ccc|ccc}
      \alpha_1 &&\alpha_2 &\beta_1 &&\beta_2
      \\ 0 && e^{2\times 2} &\beta_3&&\beta_4
      \\\hline &\gamma&&&\delta&
      \end{array}\right)=\begin{pmatrix}\alpha&\beta \\ \gamma&\delta\end{pmatrix}\]
      where $\alpha_1, \beta_2\in M_{n-2}(R_m)$, $\alpha_2, \beta_1\in M_{(n-2)\times 2}(R_m)$, $\beta_3\in M_{2}(R_m)$, $\beta_4\in M_{2\times (n-2)}(R_m)$ and $\alpha, \beta
      ,\gamma,\delta\in M_n(R_m)$. Clearly $\beta\not\equiv 0 (\textnormal{mod }I_m)$ since $\sigma\not\in U_{2n}((R_m,\Lambda_m),(I_m,I_m\cap\Lambda_m))$.\\
\\
\underline{case 2.1} Assume that $\beta_3\not\equiv 0 (\textnormal{mod }I_m)$ or $\beta_4\not\equiv 0 (\textnormal{mod }I_m)$.\\ Set $g_1:=T_{-n,n-1}(s_1)\in U_1$ and $\omega:=[\sigma^{-1},g_1]$. Then 
 \small
 \begin{align*}
 &\omega\\
 =&[\sigma^{-1},g_1]\\
 =&(e+\sigma'_{*,-n}s_1\sigma_{(n-1)*}-\sigma'_{*,-(n-1)}\lambda_ms_1\sigma_{n*})(g_1)^{-1}\\
 =&(e+s_1\bordermatrix{ & \cr 1&\bar\lambda\bar b_{n,-1} \cr &\vdots \cr n&\bar\lambda\bar b_{n,-n} \cr -n&1 \cr &0 \cr &\vdots \cr -1&0}
 \bordermatrix{&1&&&&n&-n&&-1\cr&0&\dots&0& 1 & 0 &b_{n-1,-n}& \dots & b_{n-1,-1} }\\
 &-\lambda_ms_1\bordermatrix{ & \cr 1&\bar\lambda\bar b_{n-1,-1} \cr &\vdots \cr n&\bar\lambda\bar b_{n-1,-n} \cr -n&0\cr &1 \cr &0 \cr &\vdots \cr -1&0}
 \bordermatrix{&1&&&n&-n&&-1\cr&0&\dots&0& 1&b_{n,-n}& \dots & b_{n,-1}})\\
 &\cdot (g_1)^{-1}.\hspace{9cm}\hbox{}
 \end{align*}
 \normalsize
 Since $\beta_3\not\equiv 0 (\textnormal{mod }I_m)$ or $\beta_4\not\equiv 0 (\textnormal{mod }I_m)$, $(\omega_{-n,-n},\dots,\omega_{-n,-1})\not\equiv (1,0,\dots,0)(\textnormal{mod }$ $I_m)$ or
$(\omega_{-(n-1),-n},\dots,\omega_{-(n-1),-1})\not\equiv (0,1,0,\dots,0)(\textnormal{mod }I_m)$. Hence $\omega\not\in U_{2n}((R_m,$ $\Lambda_m),(I_m,I_m\cap\Lambda_m))$. Further the next to last row of $\omega$ equals $f_{-2}$.  Set $\epsilon_2:=P_{1,-2}\in EU_{2n}(R_m,\Lambda_m)$. Then the first row of $^{\epsilon_2}\omega$ equals $f_1$. Since $U_{2n}((R_m,\Lambda_m),(I_m,I_m\cap\Lambda_m))$ is normal,  $^{\epsilon_2}\omega\not\in U_{2n}((R_m,\Lambda_m),(I_m,I_m\cap\Lambda_m))$. One can proceed now as in case 1.1 ($^{\epsilon_2}\omega$ has the same properties as $\tau$ in case 1.1).\\
\\
\underline{case 2.2} Assume that $\beta_3\equiv 0 (\textnormal{mod }I_m)$ and $\beta_4\equiv 0 (\textnormal{mod }I_m)$.\\ It follows that $\beta_1\equiv 0 (\textnormal{mod }I_m)$. Hence $\beta_2\not\equiv 0 (\textnormal{mod }I_m)$ since $\beta\not\equiv 0 (\textnormal{mod }I_m)$. Set 
\[\xi:=\prod\limits_{k=-1}^{-(n-2)}T_{n-1,k}(-b_{n-1,k})\prod\limits_{k=-1}^{-(n-2)}T_{nk}(-b_{nk})\in EU_{2n}((R_m,\Lambda_m),(I_m,\Gamma_m)).\]
      Then $\omega:=\sigma\xi$ has the form 
       \[\left(\begin{array}{ccc|ccc}
            \alpha'_1 &&\alpha'_2 &\beta'_1 &&\beta'_2
            \\ 0 && e^{2\times 2} &\beta'_3&&0
            \\\hline &\gamma'&&&\delta'&
            \end{array}\right)=\begin{pmatrix}\alpha'&\beta' \\ \gamma'&\delta'\end{pmatrix}\]
            where $\alpha'_1, \beta'_2\in M_{n-2}(R_m)$, $\alpha'_2, \beta'_1\in M_{(n-2)\times 2}(R_m)$, $\beta'_3\in M_{2}(R_m)$ and $\alpha', \beta'
            ,\gamma',\delta'\in M_n(R_m)$. Further $\alpha', \delta'\equiv e^{n\times n} (\textnormal{mod }I_m)$, $\gamma'\equiv 0 (\textnormal{mod }I_m)$, $\beta'_1\equiv 0 (\textnormal{mod }I_m)$, $\beta'_3\equiv 0 (\textnormal{mod }I_m)$ and $\beta'_2\not\equiv 0 (\textnormal{mod }I_m)$. Since $\beta'_2\not\equiv 0 (\textnormal{mod }I_m)$, there is an $i\in\{1,\dots,n-2\}$ and a $j\in\{-(n-2),\dots,-1\}$ such that $\beta'_{ij}\not\in I_m$. Choose an $l\in\{1,\dots,n-2\}\setminus\{-j\}$ and set $\epsilon_{11}:=T_{jl}(-1)\in EU_{2n}(R_m,\Lambda_m)$. Then $^{\epsilon_{11}}\omega$ has the form 
            \[\left(\begin{array}{ccc|ccc}
                  \alpha''_1 &&\alpha''_2 &\beta''_1 &&\beta''_2
                  \\ 0 && e^{2\times 2} &\beta''_3&&\beta''_4
                  \\\hline &\gamma''&&&\delta''&
                  \end{array}\right)=\begin{pmatrix}\alpha''&\beta'' \\ \gamma''&\delta''\end{pmatrix}\]
                  where $\alpha''_1, \beta''_2\in M_{n-2}(R_m)$, $\alpha''_2, \beta''_1\in M_{(n-2)\times 2}(R_m)$, $\beta''_3\in M_{2}(R_m)$, $\beta''_4\in M_{2\times (n-2)}(R_m)$ and $\alpha'', \beta''
                  ,\gamma'',\delta''\in M_n(R_m)$. Further $\alpha''_{il}\not\equiv \delta_{il}(mod~I_m)$. Hence $\alpha''\not\equiv e^{n\times n}(\textnormal{mod }I_m)$ and thus one can proceed as in case 1.1.\\
\\
\underline{Part II} Assume that $h\not\in CU_{2n}((R_m,\Lambda_m),(I_m,I_m\cap\Lambda_m))$ and $n=3$.\\
There is a $g_0\in U_0$ and $\epsilon_0,\epsilon_{1}\in EU_{6}(R_m,\Lambda_m)$ such that $\sigma:=$$^{\epsilon_{1}}[^{\epsilon_0}h,g_0]\not\in U_{6}((R_m,\Lambda_m),(I_m,I_m\cap\Lambda_m))$ and 
$\sigma$ has the form 
\[\left(\begin{array}{ccc|c}
&*& &\multirow{2}{*}{$\beta$}
\\ 0 &0& 1 &
\\\hline &\gamma&&\delta
\end{array}\right)=\begin{pmatrix}\alpha&\beta \\ \gamma&\delta\end{pmatrix}\]
where $\alpha=(\alpha_{ij})_{1\leq i,j\leq 3}, \beta=
(\beta_{ij})_{\substack{1 \leq i \leq 3\\-3 \leq j \leq -1}},\gamma=(\gamma_{ij})_{\substack{-3\leq i \leq -1\\1 \leq j \leq 3}},\delta=(\delta_{ij})_{-3 \leq i,j \leq -1}\in M_3(R_m)$ (see Part I above and \cite{bak_2}, chapter IV, Lemma 3.12, Part II, general case).\\
\\
\underline{case 1} Assume that there is an $i\in\{-3,-2,-1\}$ such that $\gamma_{i2}\not\in I_m$.\\
Set $g_1:=T_{1,-2}(s_1)\in U_1$ and $\omega:=[\sigma,g_1]$. Then 
\begin{align*}
&\omega\\
=&[\sigma,g_1]\\
=&(e+\sigma_{*1}s_1\sigma'_{-2,*}-\sigma_{*2}\bar\lambda_ms_1\sigma'_{-1,*})(g_1)^{-1}\\
=&(e+s_1\begin{pmatrix}\alpha_{11}\\\alpha_{21}\\0\\ \gamma_{-3,1}\\ \gamma_{-2,1} \\ \gamma_{-1,1}\end{pmatrix}
\begin{pmatrix}\lambda_m\bar\gamma_{-1,2}&\lambda_m\bar\gamma_{-2,2}&\lambda_m\bar\gamma_{-3,2}&0&\bar\alpha_{22}&\bar\alpha_{12}\end{pmatrix}\\
&-\bar\lambda_ms_1\begin{pmatrix}\alpha_{12}\\\alpha_{22}\\0\\ \gamma_{-3,2}\\ \gamma_{-2,2} \\ \gamma_{-1,2}\end{pmatrix}
\begin{pmatrix}\lambda_m\bar\gamma_{-1,1}&\lambda_m\bar\gamma_{-2,1}&\lambda_m\bar\gamma_{-3,1}&0&\bar\alpha_{21}&\bar\alpha_{11}\end{pmatrix})\\
&\cdot (g_1)^{-1}.
\end{align*}
\\
Assume that 
\begin{align*}&s_1\sigma_{*1}\begin{pmatrix}\lambda_m\bar\gamma_{-1,2}&\lambda_m\bar\gamma_{-2,2}&\lambda_m\bar\gamma_{-3,2}\end{pmatrix}\\
&-\bar\lambda_ms_1\sigma_{*2}\begin{pmatrix}\lambda_m\bar\gamma_{-1,1}&\lambda_m\bar\gamma_{-2,1}&\lambda_m\bar\gamma_{-3,1}\end{pmatrix}\equiv 0(\textnormal{mod }I_m).\end{align*} By multiplying $\sigma'_{1*}$ from the left we get that $s_1\begin{pmatrix}\lambda_m\bar\gamma_{-1,2}&\lambda_m\bar\gamma_{-2,2}&\lambda_m\bar\gamma_{-3,2}\end{pmatrix}$ $\equiv 0(\textnormal{mod }I_m)$ 
which implies $\begin{pmatrix}\gamma_{-1,2}&\gamma_{-2,2}&\gamma_{-3,2}\end{pmatrix}\equiv 0(\textnormal{mod }I_m)$. Since that is a contradiction, 
\begin{align*}&s_1\sigma_{*1}\begin{pmatrix}\lambda_m\bar\gamma_{-1,2}&\lambda_m\bar\gamma_{-2,2}&\lambda_m\bar\gamma_{-3,2}\end{pmatrix}\\
&-\bar\lambda_ms_1\sigma_{*2}\begin{pmatrix}\lambda_m\bar\gamma_{-1,1}&\lambda_m\bar\gamma_{-2,1}&\lambda_m\bar\gamma_{-3,1}\end{pmatrix}\not\equiv 0(\textnormal{mod }I_m)\end{align*}
and hence $\omega\not\in U_{6}((R_m,\Lambda_m),(I_m,I_m\cap\Lambda_m))$. 
Further the third row of $\omega$ equals $f_{3}$. Set $\epsilon_2:=P_{13}\in EU_{6}(R_m,\Lambda_m)$. Then the first row of $^{\epsilon_2}\omega$ equals $f_1$. Since $U_{6}((R_m,\Lambda_m),(I_m,I_m\cap\Lambda_m))$ is normal, $^{\epsilon_2}\omega\not\in U_{6}((R_m,\Lambda_m),(I_m,I_m\cap\Lambda_m))$. One can proceed now as in Part I, case 1 ($^{\epsilon_2}\omega$ has the same properties as $\tau$ in Part I, case 1).\\
\\
\underline{case 2} Assume that there is an $i\in\{-3,-2,-1\}$ such that $\gamma_{i1}\not\in I_m$.\\
This case can be treated similarly.\\
\\
\underline{case 3}
Assume that $\gamma_{-3,1},\gamma_{-3,2},\gamma_{-2,1},\gamma_{-2,2},\gamma_{-1,1},\gamma_{-1,2}\in I_m$ and one of the entries $\beta_{3,-3}$ and $\beta_{3,-2}$ does not lie in $I_m$.\\
By \cite{bak_2}, chapter IV, Lemma 3.12, Part II, general case there are $x_1, x_2\in I_m$ such that $\gamma_{-1,1}+x_2(x_1\gamma_{-1,1}+\gamma_{-2,1})\in rad (R_m)\cap I_m$ where $rad (R_m)$ is the Jacobson radical of the ring $R_m$. Set $\xi_1:=T_{-1,-2}(x_2)T_{-2,-1}(x_1)\in EU_{6}((R_m,\Lambda_m),(I_m,\Gamma_m))$. Then $\rho:=\xi_1\sigma$ has the form
\[\left(\begin{array}{ccc|c}
&*& &\multirow{2}{*}{$\beta'$}
\\ 0 &0& 1 &
\\\hline &\gamma'&&\delta'
\end{array}\right)=\begin{pmatrix}\alpha'&\beta' \\ \gamma'&\delta'\end{pmatrix}\]
where $\alpha'=(\alpha'_{ij})_{1\leq i,j\leq 3}, \beta'=(\beta'_{ij})_{\substack{1 \leq i \leq 3\\-3 \leq j \leq -1}},\gamma'=(\gamma'_{ij})_{\substack{-3\leq i \leq -1\\1 \leq j \leq 3}},\delta'=(\delta'_{ij})_{-3 \leq i,j \leq -1}\in M_3(R_m)$. Further $(\beta'_{3,-3}\not\in I_m\lor\beta'_{3,-2}\not\in I_m)\land\gamma'_{-1,1}\in rad (R_m)\cap I_m$. Set $g_1:=T_{13}(s_1)\in U_1$ and $\omega:=[\rho^{-1},g_1]$. Then
\begin{align*}
&\omega\\
=&[\rho^{-1},g_1]\\
=&(e+\rho'_{*1}s_1\rho_{3*}-\rho'_{*,-3}s_1\rho_{-1,*})(g_1)^{-1}\\
=&(e+s_1\begin{pmatrix}\bar\delta'_{-1,-1}\\\bar\delta'_{-1,-2}\\\bar\delta'_{-1,-3}\\ \lambda\bar\gamma'_{-1,3}\\ \lambda\bar\gamma'_{-1,2} \\ \lambda\bar\gamma'_{-1,1}\end{pmatrix}
\begin{pmatrix}0&0&1&\beta'_{3,-3}&\beta'_{3,-2}&\beta'_{3,-1}\end{pmatrix}\\
&-s_1\begin{pmatrix}\bar\lambda\bar\beta_{3,-1}\\\bar\lambda\bar\beta_{3,-2}\\\bar\lambda\bar\beta_{3,-3}\\1\\0\\0\end{pmatrix}
\begin{pmatrix}\gamma'_{-1,1}&\gamma'_{-1,2}&\gamma'_{-1,3}&\delta'_{-1,-3}&\delta'_{-1,-2}&\delta'_{-1,-1}\end{pmatrix})\\
&\cdot (g_1)^{-1}.
\end{align*}
Assume that 
\begin{align*}&s_1\rho'_{*1}\begin{pmatrix}\beta'_{3,-3}&\beta'_{3,-2}\end{pmatrix}\\
&-s_1\rho'_{*,-3}\begin{pmatrix}\delta'_{-1,-3}&\delta'_{-1,-2}\end{pmatrix}\equiv 0(\textnormal{mod }I_m).\end{align*} By multiplying 
$\rho_{1*}$ from the left we get that $s_1\begin{pmatrix}\beta'_{3,-3}&\beta'_{3,-2}\end{pmatrix}$ $\equiv 0(\textnormal{mod }I_m)$ 
which implies $\begin{pmatrix}\beta'_{3,-3}&\beta'_{3,-2}\end{pmatrix}\equiv 0(\textnormal{mod }I_m)$. Since that is a contradiction, 
\begin{align*}&s_1\rho'_{*1}\begin{pmatrix}\beta'_{3,-3}&\beta'_{3,-2}\end{pmatrix}\\
&-s_1\rho'_{*,-3}\begin{pmatrix}\delta'_{-1,-3}&\delta'_{-1,-2}\end{pmatrix}\not\equiv 0(\textnormal{mod }I_m)\end{align*}
and hence $\omega\not\in U_{6}((R_m,\Lambda_m),(I_m,I_m\cap\Lambda_m))$. 
Obviously $\omega_{-1,*}\equiv f_{-1}(\textnormal{mod }I_m)$ and $\omega_{-1,-1}\equiv 1 (\textnormal{mod }rad(R_m)\cap I_m)$. Set
$\epsilon_2:=P_{3,-1}\in EU_{6}(R_m,\Lambda_m)$ and $\zeta:=$$^{\epsilon_2}\omega$. Then 
$\zeta\not\in U_{6}((R_m,\Lambda_m),(I_m,I_m\cap\Lambda_m))$. Further $\zeta_{3*}\equiv f_3(\textnormal{mod }I_m)$ and $\zeta_{33}\equiv 1 (\textnormal{mod }rad(R_m)\cap I_m)$. By Nakayama's lemma $\zeta_{33}$ is invertible. Set $\xi_2:=T_{32}(-(\zeta_{33})^{-1}$ $\zeta_{32})T_{31}(-(\zeta_{33})^{-1}$ $\zeta_{31})T_{3,-1}(-(\zeta_{33})^{-1}\zeta_{3,-1})T_{3,-2}(-(\zeta_{33})^{-1}\zeta_{3,-2})\in EU_6((R_m,\Lambda_m),(I_m,$ $\Gamma_m))$ and $\eta:=\zeta\xi_2$. Then $\eta$ has the form
\[\left(\begin{array}{ccc|ccc}
&*& &&* &
\\ 0 &0& \alpha''_{33} &\beta''_{3,-3}&0&0
\\\hline &\gamma''&&&\delta''&
\end{array}\right)=\begin{pmatrix}\alpha''&\beta'' \\ \gamma''&\delta''\end{pmatrix}\]
where $\alpha''=(\alpha''_{ij})_{1\leq i,j\leq 3}, \beta''=
(\beta''_{ij})_{\substack{1 \leq i \leq 3\\-3 \leq j \leq -1}},\gamma''=(\gamma''_{ij})_{\substack{-3\leq i \leq -1\\1 \leq j \leq 3}},\delta''=(\delta''_{ij})_{-3 \leq i,j \leq -1}\in M_3(R_m)$. Since $\zeta\not\in U_{6}((R_m,\Lambda_m),(I_m,I_m\cap\Lambda_m))$, $\eta\not\in U_{6}((R_m,\Lambda_m),(I_m,I_m\cap\Lambda_m))$. Further $\alpha''_{33}\equiv 1 (\textnormal{mod }rad(R_m)\cap I_m)$ and $ \beta''_{3,-3}\in I_m$. Since $\eta_{3*}\equiv f_3(\textnormal{mod }I_m)$, $\eta_{*,-3}\equiv e_{-3}(\textnormal{mod }I_m)$ (apply Lemma 3.9 to the image of $\eta$ in $U_{2n}(R_m/I_m,$ $\Lambda_m/(\Lambda_m\cap I_m))$). Hence $\beta''_{1,-3},\beta''_{2,-3}, \delta''_{-3,-3}-1,\delta''_{-2,-3},\delta''_{-1,-3}\in I_m$.\\
\\
\underline{case 3.1} Assume that there is an $i\in\{-3,-2,-1\}$ and a $j\in\{1,2\}$ such that $\gamma''_{ij}\not\in I_m$.\\
See case 1.\\
\\
\underline{case 3.2} Assume that $\gamma''_{-3,1},\gamma''_{-3,2},\gamma''_{-2,1},\gamma''_{-2,2},\gamma''_{-1,1},\gamma''_{-1,2}\in I_m$ and one of the entries $\beta''_{1,-2}$, $\beta''_{1,-1}$, $\beta''_{2,-2}$, $\beta''_{2,-1}$, $\delta''_{-3,-2}$ and $\delta''_{-3,-1}$ does not lie in $I_m$.\\
Set $g_2:=T_{-1,2}(s_2)\in U_2$ and $\theta:=[\eta,g_2]$. Then 
\begin{align*}
&\theta\\
=&[\eta,g_2]\\
=&(e+\eta_{*,-1}s_2\eta'_{2*}-\eta_{*,-2}\lambda_m s_2\eta'_{1*})(g_2)^{-1}\\
=&(e+s_2\begin{pmatrix}\beta''_{1,-1}\\\beta''_{2,-1}\\0\\ \delta''_{-3,-1}\\ \delta''_{-2,-1} \\ \delta''_{-1,-1}\end{pmatrix}
\begin{pmatrix}\bar\delta''_{-1,-2}&\bar\delta''_{-2,-2}&\bar\delta''_{-3,-2}&0&\bar\lambda_m\bar\beta''_{2,-2}&\bar\lambda_m\bar\beta''_{1,-2}\end{pmatrix}\\
&-\lambda_ms_2\begin{pmatrix}\beta''_{1,-2}\\\beta''_{2,-2}\\0\\ \delta''_{-3,-2}\\ \delta''_{-2,-2} \\ \delta''_{-1,-2}\end{pmatrix}
\begin{pmatrix}\bar\delta''_{-1,-1}&\bar\delta''_{-2,-1}&\bar\delta''_{-3,-1}&0&\bar\lambda_m\bar\beta''_{2,-1}&\bar\lambda_m\bar\beta''_{1,-1}\end{pmatrix})\\
&\cdot (g_2)^{-1}.
\end{align*}
Assume that 
\begin{align*}&s_2\eta_{*,-1}\begin{pmatrix}\bar\delta''_{-3,-2}&0&\bar\lambda_m\bar\beta''_{2,-2}&\bar\lambda_m\bar\beta''_{1,-2}\end{pmatrix}\\
&-\lambda_m s_2\eta_{*,-2}\begin{pmatrix}\bar\delta''_{-3,-1}&0&\bar\lambda_m\bar\beta''_{2,-1}&\bar\lambda_m\bar\beta''_{1,-1}\end{pmatrix}\equiv 0(\textnormal{mod }I_m).\end{align*} It follows that  $\begin{pmatrix}\delta''_{-3,-2}&0&\beta''_{2,-2}&\beta''_{1,-2}\end{pmatrix}, \begin{pmatrix}\delta''_{-3,-1}&0&\beta''_{2,-1}&\beta''_{1,-1}\end{pmatrix} \equiv 0(\textnormal{mod }I_m)$. Since that is a contradiction, 
\begin{align*}&s_2\eta_{*,-1}\begin{pmatrix}\bar\delta''_{-3,-2}&0&\bar\lambda_m\bar\beta''_{2,-2}&\bar\lambda_m\bar\beta''_{1,-2}\end{pmatrix}\\
&-\lambda_m s_2\eta_{*,-2}\begin{pmatrix}\bar\delta''_{-3,-1}&0&\bar\lambda_m\bar\beta''_{2,-1}&\bar\lambda_m\bar\beta''_{1,-1}\end{pmatrix}\not\equiv 0(\textnormal{mod }I_m)\end{align*}
and hence $\theta\not\in U_{6}((R_m,\Lambda_m),(I_m,I_m\cap\Lambda_m))$. Clearly $\theta_{3*}=f_3$. Set $\epsilon_3:=P_{13}$. Thus the first row of $^{\epsilon_3}\theta$ equals $f_1$. Since $^{\epsilon_3}\theta\not\in U_{6}((R_m,\Lambda_m),(I_m,I_m\cap\Lambda_m))$, one can proceed now as in Part I, case 1.\\
\\
\underline{case 3.3} Assume that $\gamma''_{-3,1},\gamma''_{-3,2},\gamma''_{-2,1},\gamma''_{-2,2},\gamma''_{-1,1},\gamma''_{-1,2},\beta''_{1,-2}$, $\beta''_{1,-1}$, $\beta''_{2,-2}$, $\beta''_{2,-1}$, $\delta''_{-3,-2},\delta''_{-3,-1}\in I_m$ and one of the elements  $\alpha''_{11}-1$, $\alpha''_{21}$, $\gamma''_{-1,3}$, $\delta''_{-1,-1}-1$ and $\delta''_{-1,-2}$ does not lie in $I_m$.\\
Set $g_2:=T_{13}(s_2)\in U_2$ and $\theta:=[\eta^{-1},g_2]$. Then 
\begin{align*}
&\theta\\
=&[\eta^{-1},g_2]\\
=&(e+\eta'_{*1}s_2\eta_{3*}-\eta'_{*,-3}s_2\eta_{-1,*})(g_2)^{-1}\\
=&(e+s_2\begin{pmatrix}\bar\delta''_{-1,-1}\\\bar\delta''_{-1,-2}\\\bar\delta''_{-1,-3}\\\lambda_m\bar\gamma''_{-1,3}\\ \lambda_m\bar\gamma''_{-1,2} \\ \lambda_m\bar\gamma''_{-1,1}\end{pmatrix}
\begin{pmatrix}0&0&\alpha''_{33}&\beta''_{3,-3}&0&0\end{pmatrix}\\
&-s_2\begin{pmatrix}0\\0\\ \bar\lambda_m\bar\beta''_{3,-3}\\ \bar\alpha''_{33}\\0\\0\end{pmatrix}
\begin{pmatrix}\gamma''_{-1,1}&\gamma''_{-1,2}&\gamma''_{-1,3}&\delta''_{-1,-3}&\delta''_{-1,-2}&\delta''_{-1,-1}\end{pmatrix})\\
&\cdot (g_2)^{-1}.
\end{align*}
Assume that $\theta\in U_{6}((R_m,\Lambda_m),(I_m,I_m\cap\Lambda_m))$. Then $s_2\bar\delta''_{-1,-1}\alpha''_{33}-s_2=\theta_{13}\in I_m$ and  $s_2\bar\delta''_{-1,-2}\alpha''_{33}=\theta_{23}\in I_m$. 
Since $\alpha''_{33}\equiv 1 (\textnormal{mod }I_m)$, it follows that $\delta''_{-1,-1}-1, \delta''_{-1,-2}\in I_m$. Consider the column \[\eta'_{*1}s_2\alpha''_{33}-\eta'_{*,-3}s_2\gamma''_{-1,3}-\begin{pmatrix}s_2&0&-s_2^2\bar\lambda_m\bar\beta''_{3,-3}\gamma''_{-1,1}&-s_2^2\bar\alpha''_{33}\gamma''_{-1,1} &0&0\end{pmatrix}^t.\]
Since by assumption $\theta\in U_{6}((R_m,\Lambda_m),(I_m,I_m\cap\Lambda_m))$, the column above is congruent to $0$ modulo $I_m$. By multiplying $\eta_{-3,*}$ from the left we get that $\gamma''_{-1,3}\in I_m$ since $\gamma''_{-3,1}, \gamma''_{-1,1}\in I_m$ ($\gamma''_{-3,1}$ is the first entry of $\eta_{-3,*}$). Hence $\eta_{-1,*}\equiv f_{-1}(\textnormal{mod }I_m)$. It follows that $\eta_{*1}\equiv e_1(\textnormal{mod }I_m)$ (apply Lemma 3.9 to the image of $\eta$ in $U_{2n}(R_m/I_m,$ $\Lambda_m/(\Lambda_m\cap I_m))$) and hence $\alpha''_{11}-1,\alpha''_{21}\in I_m$. Since that is a contradiction, $\theta\not\in U_{6}((R_m,\Lambda_m),(I_m,I_m\cap\Lambda_m))$.\\
\\
\underline{case 3.3.1} Assume that $\theta_{13}\not\in I_m$ or $\theta_{23}\not\in I_m$.\\
Set $g_3:=T_{-2,1}(s_3)\in U_3$ and $\tau:=[\theta^{-1},g_3]$. Then 
\begin{align*}
&\tau\\
=&[\theta^{-1},g_3]\\
=&(e+\theta'_{*,-2}s_3\theta_{1*}-\theta'_{*,-1}\lambda_ms_3\theta_{2*})(g_3)^{-1}
\end{align*}
\begin{align*}
=&(e+s_3\begin{pmatrix}\bar\lambda_m\bar\theta_{2,-1}\\0\\\bar\lambda_m\bar\theta_{2,-3}\\\bar\theta_{23}\\\bar\theta_{22}\\\bar\theta_{21}\end{pmatrix}
\begin{pmatrix}\theta_{11}&\theta_{12}&\theta_{13}&\theta_{1,-3}&0&\theta_{1,-1}\end{pmatrix}\\
&-\lambda_ms_3\begin{pmatrix}\bar\lambda_m\bar\theta_{1,-1}\\0\\\bar\lambda_m\bar\theta_{1,-3}\\\bar\theta_{13}\\\bar\theta_{12}\\\bar\theta_{11}\end{pmatrix}
\begin{pmatrix}\theta_{21}&\theta_{22}&\theta_{23}&\theta_{2,-3}&0&\theta_{2,-1}\end{pmatrix})\\
&\cdot (g_3)^{-1}.
\end{align*}
Assume that $\tau\in U_{6}((R_m,\Lambda_m),(I_m,I_m\cap\Lambda_m))$. Then $\theta'_{*,-2}s_3\theta_{13}-\theta'_{*,-1}\lambda_ms_3\theta_{23}\equiv 0 (\textnormal{mod }I_m)$. 
It follows that $\theta_{13},\theta_{23}\in I_m$ which is a contradiction. Hence $\tau\not\in U_{6}((R_m,\Lambda_m),(I_m,I_m\cap\Lambda_m))$. Clearly $\tau_{2*}=f_2$. Set $\epsilon_4:=P_{12}$. Then the first row of $^{\epsilon_4}\theta$ equals $f_1$. Since $^{\epsilon_4}\theta\not\in U_{6}((R_m,\Lambda_m),(I_m,I_m\cap\Lambda_m))$, one can proceed now as in Part I, case 1.\\
\\
\underline{case 3.3.2} Assume that $\theta_{13}\in I_m$ and $\theta_{23}\in I_m$.\\
Let $\hat\theta$ be the image of $\theta$ in $U_{2n}(R_m/I_m,\Lambda_m/(\Lambda_m\cap I_m))$. Clearly $\hat\theta$ has the form
{\renewcommand{\arraystretch}{1.3}
\[\left(\begin{array}{ccc|ccc}
1&0&0 &0&0&0
\\ 0 &1& 0 &0&0&0
\\ 0&0&\hat\theta_{33}&0&0&0
\\\hline 0&0&\hat\theta_{-3,3}&\hat\theta_{-3,-3}&\hat\theta_{-3,-2}&\hat\theta_{-3,-1}
\\ 0&0&0&0&1&0
\\ 0&0&0&0&0&1
\end{array}\right).\]
}
It follows that $\hat\theta_{33}$ is invertible. Let $\mathbbm{h}$ be the map defined in Definition 3.3. Then
\begin{align*}
&\bar{\hat\theta}_{33}\hat\theta_{-3,-2}\\
=&\mathbbm{h}(\hat\theta_{*3},\hat\theta_{*,-2})\\
=&\mathbbm{h}(\hat\theta e_{3},\hat\theta e_{-2})\\
=&\mathbbm{h}(e_{3},e_{-2})\\
=&0.
\end{align*}
Hence $\hat\theta_{-3,-2}=0$ and therefore $\theta_{-3,-2}\in I_m$. Further 
\begin{align*}
&\bar{\hat\theta}_{33}\hat\theta_{-3,-1}\\
=&\mathbbm{h}(\hat\theta_{*3},\hat\theta_{*,-1})\\
=&\mathbbm{h}(\hat\theta e_{3},\hat\theta e_{-1})\\
=&\mathbbm{h}(e_{3},e_{-1})\\
=&0.
\end{align*}
Hence $\hat\theta_{-3,-1}=0$ and therefore $\theta_{-3,-1}\in I_m$. Clearly $\theta_{22}=1$. Set $\xi_3:=T_{23}(-\theta_{23})$ $T_{2,-1}(-\theta_{2,-1})\in U_{6}((R_m,\Lambda_m),(I_m,\Gamma_m))$ and $\chi:=\theta\xi_3$. Then $\chi_{23}, \chi_{2,-1}=0$ 
and the image $\hat\chi$ of $\chi$ in $U_{6}(R_m/I_m,\Lambda_m/(\Lambda_m\cap I_m))$ has the form 
{\renewcommand{\arraystretch}{1.3}
\[\left(\begin{array}{ccc|ccc}
1&0&0 &0&0&0
\\ 0 &1& 0 &0&0&0
\\ 0&0&\hat\chi_{33}&0&0&0
\\\hline 0&0&\hat\chi_{-3,3}&\hat\chi_{-3,-3}&0&0
\\ 0&0&0&0&1&0
\\ 0&0&0&0&0&1
\end{array}\right).\]
}
Since $\theta\not\in U_{6}((R_m,\Lambda_m),(I_m,I_m\cap\Lambda_m))$, $\chi\not\in U_{6}((R_m,\Lambda_m),(I_m,I_m\cap\Lambda_m))$ (i.e. $\hat\chi\neq e$). Set $g_3:=T_{31}(s_3)\in U_3$ and $\mu:=[\chi,g_3]$. Then 
\begin{align*}
&\mu\\
=&[\chi,g_3]\\
=&(e+\chi_{*3}s_3\chi'_{1*}-\chi_{*,-1}s_3\chi'_{-3,*})(g_3)^{-1}\\
=&(e+s_3\begin{pmatrix}\chi_{13}\\0\\\chi_{33}\\\chi_{-3,3}\\\chi_{-2,3}\\\chi_{-1,3}\end{pmatrix}
\begin{pmatrix}\bar\chi_{-1,-1}&\bar\chi_{-2,-1}&\bar\chi_{-3,-1}&\bar\lambda_m\bar\chi_{3,-1}&0&\bar\lambda_m\bar\chi_{1,-1}\end{pmatrix}\\
&-s_3\begin{pmatrix}\chi_{1,-1}\\0\\\chi_{3,-1}\\\chi_{-3,-1}\\\chi_{-2,-1}\\\chi_{-1,-1}\end{pmatrix}
\begin{pmatrix}\lambda_m\bar\chi_{-1,3}&\lambda_m\bar\chi_{-2,3}&\lambda_m\bar\chi_{-3,3}&\bar\chi_{33}&0&\bar\chi_{13}\end{pmatrix})\\
&\cdot (g_3)^{-1}.
\end{align*}
Assume that $\mu\in U_{6}((R_m,\Lambda_m),(I_m,I_m\cap\Lambda_m))$. Then $s_3\chi_{33}\bar\chi_{-1,-1}-s_3\chi_{3,-1}\lambda_m\bar\chi_{-1,3}-s_3(1+s_3\chi_{33}\bar\chi_{-3,-1}-s_3\chi_{3,-1}\lambda_m\bar\chi_{-3,3})=\mu_{31}\in I_m$ and hence $s_3\chi_{33}\bar\chi_{-1,-1}-s_3\in I_m$. It follows that $\chi_{33}\equiv 1(\textnormal{mod }I_m)$ (i.e. $\hat\chi_{33}=1$) since $s_3\in (R_m)^*$ and $\chi_{-1,-1}\equiv 1 (\textnormal{mod }I_m)$. That implies
\begin{align*}
&\hat\chi_{-3,-3}\\
=&\mathbbm{h}(\hat\chi_{*3},\hat\chi_{*,-3})\\
=&\mathbbm{h}(\hat\chi e_{3},\hat\chi e_{-3})\\
=&\mathbbm{h}(e_{3},e_{-3})\\
=&1.
\end{align*}
Further $s_3\chi_{-3,3}\bar\chi_{-1,-1}-s_3\chi_{-3,-1}\lambda_m\bar\chi_{-1,3}-s_3(s_3\chi_{-3,3}\bar\chi_{-3,-1}-s_3\chi_{-3,-1}\lambda_m\bar\chi_{-3,3})=\mu_{-3,1}\in I_m$ and hence $s_3\chi_{-3,3}\bar\chi_{-1,-1}\in I_m$. It follows that $\chi_{-3,3}\equiv 0(\textnormal{mod }I_m)$ (i.e. $\hat\chi_{-3,3}=0$) since $s_3\in (R_m)^*$ and $\chi_{-1,-1}\equiv 1 (\textnormal{mod }I_m)$. But that implies the contradiction $\hat\chi=e$. Hence $\mu\not\in U_{6}((R_m,\Lambda_m),(I_m,I_m\cap\Lambda_m))$. Clearly $\mu_{2*}=f_2$. Set $\epsilon_4:=P_{12}$. Then the first row of $^{\epsilon_4}\mu$ equals $f_1$. Since $^{\epsilon_4}\mu\not\in U_{6}((R_m,\Lambda_m),(I_m,I_m\cap\Lambda_m))$, one can proceed as in Part I, case 1.\\
\\
\underline{case 3.4} Assume that $\gamma''_{-3,1},\gamma''_{-3,2},\gamma''_{-2,1},\gamma''_{-2,2},\gamma''_{-1,1},\gamma''_{-1,2},\beta''_{1,-2}$, $\beta''_{1,-1}$, $\beta''_{2,-2}$, $\beta''_{2,-1}$, $\delta''_{-3,-2},\delta''_{-3,-1},\alpha''_{11}-1$, $\alpha''_{21}$, $\gamma''_{-1,3}$, $\delta''_{-1,-1}-1,\delta''_{-1,-2}\in I_m$ and one of the elements $\alpha''_{12}, \alpha''_{22}-1, \gamma''_{-2,3}, \delta''_{-2,-1}, \delta''_{-2,-2}-1, \alpha''_{13}$ and $\alpha''_{23}$ does not lie in $I_m$.\\
Set $g_2:=T_{21}(s_2)\in U_2$ and $\theta:=[\eta,g_2]$. Then 
\begin{align*}
&\theta\\
=&[\eta,g_2]\\
=&(e+\eta_{*2}s_2\eta'_{1*}-\eta_{*,-1}s_2\eta'_{-2,*})(g_2)^{-1}\\
=&(e+s_2\begin{pmatrix}\alpha''_{12}\\\alpha''_{22}\\0\\ \gamma''_{-3,2}\\ \gamma''_{-2,2} \\ \gamma''_{-1,2}\end{pmatrix}
\begin{pmatrix}\bar\delta''_{-1,-1}&\bar\delta''_{-2,-1}&\bar\delta''_{-3,-1}&0&\bar\lambda_m\bar\beta''_{2,-1}&\bar\lambda_m\bar\beta''_{1,-1}\end{pmatrix}\\
&-s_2\begin{pmatrix}\beta''_{1,-1}\\\beta''_{2,-1}\\0\\ \delta''_{-3,-1}\\ \delta''_{-2,-1} \\ \delta''_{-1,-1}\end{pmatrix}
\begin{pmatrix}\lambda_m\bar\gamma''_{-1,2}&\lambda_m\bar\gamma''_{-2,2}&\lambda_m\bar\gamma''_{-3,2}&0&\bar\alpha''_{22}&\bar\alpha''_{12}\end{pmatrix})\\
&\cdot (g_2)^{-1}.
\end{align*}
Assume that $\theta\in U_{6}((R_m,\Lambda_m),(I_m,I_m\cap\Lambda_m))$. Then 
$s_2\eta_{*2}\bar\lambda_m\bar\beta''_{1,-1}-s_2\eta_{*,-1}\bar\alpha''_{12}\equiv 0(\textnormal{mod }I_m)$. It follows that $\alpha''_{12}\in I_m$. Hence $-s_2\eta_{*,-1}\bar\alpha''_{22}+s_2e_{-1}\equiv 0(\textnormal{mod }I_m)$. By multiplying 
$\eta'_{-1,*}$ from the left we get that $-s_2\bar\alpha''_{22}+s_2\bar\alpha''_{11}\in I_m$ which implies $\alpha''_{22}\equiv 1 (\textnormal{mod }I_m)$ since $\alpha''_{11}\equiv 1(\textnormal{mod }I_m)$. Let $\hat\eta$ be the image of $\eta$ in $U_{2n}(R_m/I_m,\Lambda_m/(\Lambda_m\cap I_m))$. By Lemma 3.9, $\hat\eta_{-2,*}=f_{-2}$ since $\hat\eta_{*2}=e_2$. Hence $\gamma''_{-2,3},\delta''_{-2,-2}-1,\delta''_{-2,-1}\in I_m$. Hence $\hat\eta$ has the form 
{\renewcommand{\arraystretch}{1.3}
\[\left(\begin{array}{ccc|ccc}
1&0&\hat\alpha''_{13} &0&0&0
\\ 0 &1& \hat\alpha''_{23} &0&0&0
\\ 0&0&1&0&0&0
\\\hline 0&0&\hat\gamma''_{-3,3}&1&0&0
\\ 0&0&0&0&1&0
\\ 0&0&0&0&0&1
\end{array}\right).\]
Clearly
\begin{align*}
&\overline{\hat\alpha''_{13}}\\
=&\mathbbm{h}(\hat\eta_{*3},\hat\eta_{*,-1})\\
=&\mathbbm{h}(\hat\eta e_{3},\hat\eta e_{-1})\\
=&\mathbbm{h}(e_{3},e_{-1})\\
=&0
\end{align*}
and
\begin{align*}
&\overline{\hat\alpha''_{23}}\\
=&\mathbbm{h}(\hat\eta_{*3},\hat\eta_{*,-2})\\
=&\mathbbm{h}(\hat\eta e_{3},\hat\eta e_{-2})\\
=&\mathbbm{h}(e_{3},e_{-2})\\
=&0.
\end{align*}
Hence $\hat\alpha''_{13}=0=\hat\alpha''_{23}$ and therefore $\alpha''_{13},\alpha''_{23}\in I_m$.
Since that is a contradiction, $\theta\not\in U_{6}((R_m,\Lambda_m),(I_m,I_m\cap\Lambda_m))$. Clearly $\theta_{3*}=f_3$. Set $\epsilon_3:=P_{13}$. Then the first row of $^{\epsilon_3}\theta$ equals $f_1$. Since $^{\epsilon_3}\theta\not\in U_{6}((R_m,\Lambda_m),(I_m,I_m\cap\Lambda_m))$, one can proceed as in Part I, case 1.\\
\\
\underline{case 3.5} Assume that $\gamma''_{-3,1},\gamma''_{-3,2},\gamma''_{-2,1},\gamma''_{-2,2},\gamma''_{-1,1},\gamma''_{-1,2},\beta''_{1,-2}$, $\beta''_{1,-1}$, $\beta''_{2,-2}$, $\beta''_{2,-1}$, $\delta''_{-3,-2},\delta''_{-3,-1},\alpha''_{11}-1$, $\alpha''_{21}$, $\gamma''_{-1,3}$, $\delta''_{-1,-1}-1,\delta''_{-1,-2}, \alpha''_{12}, \alpha''_{22}-1, \gamma''_{-2,3}, \delta''_{-2,-1}, \delta''_{-2,-2}-1, \alpha''_{13},\alpha''_{23}\in I_m$.\\
Since $\eta\not\in U_{6}((R_m,\Lambda_m),(I_m,I_m\cap\Lambda_m))$, $\gamma''_{-3,3}\not\in I_m$. By \cite{bak_2}, chapter IV, Lemma 3.12, Part II, case 4 there are $x_1, x_2\in I_m$ such that $\gamma''_{-1,3}+x_2(x_1\gamma''_{-1,3}+\gamma''_{-2,3})\in rad (R_m)\cap I_m$. Set $\xi_3:=T_{-1,-2}(x_2)T_{-2,-1}(x_1)\in EU_{6}((R_m,\Lambda_m),(I_m,\Gamma_m))$ and $\theta:=\xi_3\eta$. Then $\theta\equiv\eta (\textnormal{mod }I_m)$, $\theta_{33}=\alpha''_{33}\equiv 1 (\textnormal{mod }rad (R_m)\cap I_m)$ and $\theta_{-1,3}\in rad (R_m)\cap I_m$. Set $\epsilon_{21}:=T_{3,-1}(-1)\in EU_6(R_m,\Lambda_m)$ and $\tau:=$$^{\epsilon_{21}}\theta$. Then $\tau_{33}\equiv 1 (\textnormal{mod }rad (R_m)\cap I_m)$ and $\tau_{-3,-1}=\theta_{-3,-1}-\theta_{-3,3}\equiv \theta_{-3,3}\equiv\gamma''_{-3,3}(\textnormal{mod }I_m)$. Since $\gamma''_{-3,3}\not\in I_m$, it follows that $\tau_{-3,-1}\not\in I_m$. Set $\xi_4:=T_{32}(-(\tau_{33})^{-1}\tau_{32})T_{31}(-(\tau_{33})^{-1}$ $\tau_{31})T_{3,-1}(-(\tau_{33})^{-1}\tau_{3,-1})T_{3,-2}(-(\tau_{33})^{-1}\tau_{3,-2})\in EU_6((R_m,\Lambda_m), (I_m,\Gamma_m))$ and $\chi:=\tau\xi_4$. Then $\chi$ has the form
\[\left(\begin{array}{ccc|ccc}
&*& &&* &
\\ 0 &0& \alpha'''_{33} &\beta'''_{3,-3}&0&0
\\\hline &\gamma'''&&&\delta'''&
\end{array}\right)=\begin{pmatrix}\alpha'''&\beta''' \\ \gamma'''&\delta'''\end{pmatrix}\]
where $\alpha''', \beta''',\gamma''',\delta'''\in M_3(R_m)$. Further $\alpha'''_{33}\equiv 1 (\textnormal{mod }rad(R_m)\cap I_m)$, $\beta'''_{3,-3}\in I_m$ and $\delta'''_{-3,-1}\not\in I_m$. One can proceed now as in case 3.1 or case 3.2.\\
\\
\underline{case 4}
Assume that $\gamma_{-3,1},\gamma_{-3,2},\gamma_{-2,1},\gamma_{-2,2},\gamma_{-1,1},\gamma_{-1,2},\beta_{3,-3},\beta_{3,-2}\in I_m$ and $\beta_{3,-1}\not\in I_m$.\\
Set $\epsilon_{11}:=T_{21}(1)\in EU_6(R_m,\Lambda_m)$ and $\rho:=$$^{\epsilon_{11}}\sigma$. Clearly $\rho_{-3,1},\rho_{-3,2},\rho_{-2,1},\rho_{-2,2},\rho_{-1,1},$ $\rho_{-1,2}\in I_m$. Further \[\rho_{3*}=\begin{pmatrix}0&0&1&\beta_{3,-3}&\beta_{3,-2}+\beta_{3,-1}&\beta_{3,-1}\end{pmatrix}.\] Since $\beta_{3,-2}\in I_m$ and $\beta_{3,-1}\not\in I_m$, $\beta_{3,-2}+\beta_{3,-1}\not\in I_m$. One can proceed now as in case 3.\\
\\
\underline{case 5} Assume that $\gamma_{-3,1},\gamma_{-3,2},\gamma_{-2,1},\gamma_{-2,2},\gamma_{-1,1},\gamma_{-1,2},\beta_{3,-3},\beta_{3,-2},\beta_{3,-1}\in I_m$.\\
One can proceed as in case 3 ($\sigma$ has the same properties as $\zeta$ in case 3).\\
\\
\underline{Part III} Assume that $h\in CU_{2n}((R_m,\Lambda_m),(I_m,I_m\cap\Lambda_m))$.\\
This part corresponds to Proposition 3.3 in \cite{bak_2}, chapter IV. By \cite{bak_2}, chapter IV, Corollary 3.4, there is an elementary short or long root element $T_{ij}(x)\in EU_{2n}(R_m,\Lambda_m)$ such that $[h,T_{ij}(x)]\not\in U_{2n}((R_m,\Lambda_m),(I_m,\Gamma_m))$ since $h\not\in CU_{2n}((R_m,\Lambda_m),(I_m,\Gamma_m))$. Since $h\in CU_{2n}((R_m,\Lambda_m),$ $(I_m,I_m\cap\Lambda_m))$, $[h,T_{ij}(x)]\in U_{2n}((R_m,\Lambda_m),(I_m,I_m\cap\Lambda_m))$. It follows from Lemma 4.6 that there is an elementary matrix $g_0\in U_0$ such that $[h,g_0]\in U_{2n}((R_m,\Lambda_m),(I_m,I_m\cap\Lambda_m))\setminus U_{2n}((R_m,\Lambda_m),(I_m,\Gamma_m))$ (note that $J(h)\subseteq I_m$ and $h_{kk}\equiv h_{ll}(mod~I_m)~\forall k,l$ since $h\in CU_{2n}((R_m,\Lambda_m),(I_m,I_m\cap\Lambda_m))$). Set $\sigma:=[h,g_0]$. By Lemma 4.7, there is an $\epsilon_1\in EU_{2n}(R_m,\Lambda_m)$ such that $y_1:=(^{\epsilon_1}\sigma)_{11}$ is invertible. Clearly $^{\epsilon_1}\sigma\in U_{2n}((R_m,\Lambda_m),(I_m,I_m\cap\Lambda_m))\setminus U_{2n}((R_m,\Lambda_m),(I_m,\Gamma_m))$. Set $\omega:={}^{\epsilon_1}\sigma$, $\xi_1:=
T_{-2,1}(-\omega_{-2,1}(y_1)^{-1})\dots T_{21}(-\omega_{21}(y_1)^{-1})\in EU_{2n}((R_m,\Lambda_m),(I_m,\Gamma_m))$ and $\xi_2:=
T_{12}(-(y_1)^{-1}\omega_{12})\dots T_{1,-2}(-(y_1)^{-1}\omega_{1,-2})\in EU_{2n}((R_m,\Lambda_m),(I_m,\Gamma_m))$.Then $\tau:=\xi_1\omega\xi_2$ has the form 
\[\begin{pmatrix}y_1&0&y_2\\0&A&0\\y_3&0&y_4\end{pmatrix}\] where $y_2,y_3,y_4\in R_m$ and $A\in M_{2n-2}(R_m)$. \\
\\
\underline{case 1} Assume that $y_3(y_1)^{-1}\in \Gamma_m$ and $(y_1)^{-1}y_2\in \bar\lambda_m\Gamma_m$.\\ Set $\xi_3:=T_{-1,1}(-y_3(y_1)^{-1})\in EU_{2n}((R_m,\Lambda_m),(I_m,\Gamma_m))$ and $\xi_4:=T_{1,-1}(-(y_1)^{-1}y_2)\in EU_{2n}((R_m,\Lambda_m),(I_m,\Gamma_m))$. Then $\zeta:=\xi_3\tau\xi_4$ has the form \[\begin{pmatrix}y_1&0&0\\0&A&0\\0&0&y_5\end{pmatrix}\] where $y_5\in R_m$. Clearly $\zeta\in U_{2n}((R_m,$ $\Lambda_m),(I_m,I_m\cap\Lambda_m))$ but $\zeta\not\in U_{2n}((R_m,$ $\Lambda_m),(I_m,\Gamma_m))$ . Hence there is an $l\in\{2,\dots,$ $-2\}$ such that $|\zeta_{*l}|\not \in\Gamma_m$.  \\
\\
\underline{case 1.1} Assume that $\epsilon(l)=1$.\\ There are a $b'\in R$ and a $t'\in S_m$ such that $y_1|\zeta_{*l}|\bar y_1-\zeta_{-l,l}\bar y_1+\lambda_m\overline{\zeta_{-l,l}\bar y_1}=\frac{b'}{t'}$. Set $t:=\frac{t'}{1}\in R_m$ and $g_1:=T_{l,-1}(s_1s_2t)\in U_1$. One can show that $[\zeta,g_1]$ equals
\begin{align*}
&T_{2,-1}(s_1s_2t\zeta_{2l}\bar y_1)\cdot\\
&\vdots\\
&T_{l-1,-1}(s_1s_2t\zeta_{(l-1)l}\bar y_1)\cdot\\
&T_{l,-1}(s_1s_2t\zeta_{ll}\bar y_1-s_1s_2t)\cdot\\
&T_{l+1,-1}(s_1s_2t\zeta_{(l+1)l}\bar y_1)\cdot\\
&\vdots\\
&T_{-2,-1}(s_1s_2t\zeta_{-2,l}\bar y_1)\cdot\\
&T_{1,-1}(z)\\
\end{align*}
where $z=\bar\lambda_m\overline{s_1s_2t}(y_1|\zeta_{*l}|\bar y_1-\zeta_{-l,l}\bar y_1+\lambda_m\overline{\zeta_{-l,l}\bar y_1})s_1s_2t$. Since $|\zeta_{*l}|\not\in\Gamma_m$ and $y_1$ is invertible, $y_1|\zeta_{*l}|\bar y_1\not\in\Gamma_m$. Since $\zeta_{-l,l}\in I_m$, $-\zeta_{-l,l}\bar y_1+\lambda_m\overline{\zeta_{-l,l}\bar y_1}\in(\Gamma_m)_{min}\subseteq\Gamma_m$. Since $s_1s_2t$ is invertible, it follows that $\overline{s_1s_2t}(y_1|\zeta_{*l}|\bar y_1-\zeta_{-l,l}\bar y_1+\lambda_m\overline{\zeta_{-l,l}\bar y_1})s_1s_2t\not\in\Gamma_m$. Hence $z\not\in\bar\lambda_m\Gamma_m$. Set 
\begin{align*}
&\xi_5:=\\
&T_{-2,-1}(-s_1s_2t\zeta_{-2,l}\bar y_1)\cdot\\
&\vdots\\
&T_{l+1,-1}(-s_1s_2t\zeta_{(l+1)l}\bar y_1)\cdot\\
&T_{l,-1}(-(s_1s_2t\zeta_{ll}\bar y_1-s_1s_2t))\cdot\\
&T_{l-1,-1}(-s_1s_2t\zeta_{(l-1)l}\bar y_1)\cdot\\
&\vdots\\
&T_{2,-1}(-s_1s_2t\zeta_{2l}\bar y_1)\in EU_{2n}((R_m,\Lambda_m),(I_m,\Gamma_m))
\end{align*}
and $g_2:=T_{1,-1}(z)\in U_2$. Then  
\[(\xi_5[\xi_3\xi_1(^{\epsilon_1}[h,g_0])\xi_2\xi_4,g_1])=g_2.\]
Note that $g_2\not\in U_{2n}((R_m,\Lambda_m),(I_m,\Gamma_m))$ since $z\not\in\bar\lambda_m\Gamma_m$.
Set $g'_i:=\psi_m(g_i)~\forall i\in\{0,1,2\}$ and $\epsilon'_1:=\psi_m(\epsilon_1)$. Then 
\[[^{\epsilon'_1}[h',g'_0],g'_1]=g'_2.\]\\
\underline{case 1.2} Assume that $\epsilon(l)=-1$.\\ This case can be treated similarly.\\
\\
\underline{case 2} Assume that $y_3(y_1)^{-1}\not\in \Gamma_m$. \\
\\
\underline{case 2.1} Assume that $|\tau_{*l}|\in\Gamma_m~\forall l\in\{2,\dots,-2\}$.\\ Set $\chi:=T_{-1,1}(-y_3(y_1)^{-1})\in EU_{2n}(R_m,$ $\Lambda_m)$ (one checks easily that $-y_3(y_1)^{-1}\in\Lambda_m$). Then $\zeta:=\chi\tau$ has the form \[\begin{pmatrix}y_1&0&y_2\\0&A&0\\0&0&y_5\end{pmatrix}\] where $y_5\in R_m$. There is an $b'\in R$ and a $t'\in S_m$ such that $y_3(y_1)^{-1}=\frac{b'}{t'}$. Set $t:=\frac{t'}{1}\in R_m$ and $g_1:=T_{12}(s_1s_2t)\in U_1$. Using the equality $[\alpha\beta,\gamma]=$ $^{\alpha}[\beta,\gamma][\alpha,\gamma]$ one gets that $[\tau,g_1]=[\chi^{-1}\zeta,g_1]=$ $^{\chi^{-1}}[\zeta,g_1][\chi^{-1},g_1]$. It is easy to show that $[\zeta,g_1]\in EU_{2n}((R_m,\Lambda_m),(I_m,\Gamma_m))$ and hence $^{\chi^{-1}}[\zeta,g_1]\in EU_{2n}((R_m,\Lambda_m),(I_m,\Gamma_m))$. On the other hand $[\chi^{-1},g_1]=T_{-1,2}(y_3(y_1)^{-1}s_1s_2t)T_{-2,2}(-\overline{s_1s_2t}y_3(y_1)^{-1}s_1s_2t)$, by $(R6)$ in Lemma 3.12. Set $\xi_3:=T_{-1,2}(-y_3(y_1)^{-1}s_1s_2t)(^{\chi^{-1}}[\zeta,g_1])^{-1}\in EU_{2n}((
R_
m,\Lambda_m),(I_m,$ $\Gamma_m))$ and $g_2:=T_{-2,2}($ $-\overline{s_1s_2t}y_3$ $(y_1)^{-1}s_1s_2t)\in U_2$. Since 
$y_3(y_1)^{-1}\not\in \Gamma_m$, $g_2\not\in U_{2n}((R_m,\Lambda_m),(I_m,\Gamma_m))$. Clearly 
\[(\xi_3[\xi_1(^{\epsilon_1}[h,g_0])\xi_2,g_1])=g_2.\] 
As above, push this equation into $U_{2n}(R_m,\Lambda_m)/U_{2n}((R_m,\Lambda_m),(I_m,\Gamma_m))$ by applying $\psi_m$.\\
\\
\underline{case 2.2} Assume that there is an $l\in\{2,\dots,-2\}$ such that $|\tau_{*l}|\not\in\Gamma_m$.\\ Choose a $p\in\{2,\dots,-2\}$ such that $p\neq \pm l$ and set $g_1:=T_{lp}(s_1)\in U_1$. Then $\zeta:=[\tau,g_1]$ has the form \[\begin{pmatrix} 1&0&0\\ 0&B&0 \\ 0&0&1\end{pmatrix}\] where $B\in M_{2n-2}(R_m)$. Since $\tau\in U_{2n}((R_m,\Lambda_m),(I_m,I_m\cap\Lambda_m))$, $\zeta\in U_{2n}((R_m,\Lambda_m),$ $(I_m,I_m\cap\Lambda_m))$. By Lemma 4.6, $|\zeta_{*p}|\not\in\Gamma_m$. Hence $\zeta\not\in U_{2n}((R_m,\Lambda_m),(I_m,\Gamma_m))$ and thus one can proceed as in case 1.\\
\\
\underline{case 3} Assume that $(y_1)^{-1}y_2\not\in \bar\lambda_m\Gamma_m$.\\ See case 2.
\hbox{}\hfill $\qed$}
\Lemma{Let $(I,\Omega)$ be a form ideal of $(R,\Lambda)$ and $m$ a maximal ideal of $C$ such that $I\cap C\subseteq m$. Then the following is true.
\begin{enumerate}[(4.9.1)]
\item If $U' \in A'$ and $g'\in \phi_m(\psi(EU_{2n}(R,\Lambda)))$ is the nontrivial image of an elementary matrix in $EU_{2n}(R,\Lambda)$, then 
\[V'\subseteq{}^{U'}g'\]
for some $V'\in B'$.
\item If $V'\in B'$ and $d'\in \psi_m(EU_{2n}(R_m,\Lambda_m))$ is the  image of an elementary matrix in $EU_{2n}(R_m,\Lambda_m)$, then 
\[g'\in{}^{d'}V'\]
for some nontrivial image $g'\in \phi_m(\psi(EU_{2n}(R,\Lambda)))$ of an elementary matrix in $EU_{2n}(R,\Lambda)$.
\end{enumerate}
}
\Proof{Follows from the relations (R1)-(R6) in Lemma 3.12.\hfill $\qed$}
\Corollary{Let $(I,\Omega)$ be a form ideal of $(R,\Lambda)$ and $m$ a maximal ideal of $C$ such that $I\cap C\subseteq m$. If $U' \in A'$, $d'\in \psi_m(EU_{2n}(R_m,\Lambda_m))$ and $g'\in \phi_m(\psi(EU_{2n}(R,$ $\Lambda)))$ is the nontrivial image of an elementary matrix in $EU_{2n}(R,\Lambda)$, then 
\[V'\subseteq{}^{U'}(^{d'}g')\]
for some $V'\in B'$.}
\Proof{If $d' = 1$, then we are done, by  $(4.9.1)$. Assume $d' \neq 1$ and write $d'$ as a product $d_k'\dots d_1'$ of nontrivial images of elementary matrices in $EU_{2n}(R_m,\Lambda_m)$. We proceed by induction on $k$.\\
\underline{case 1} Assume that $k = 1$. Since $(A',B')$ is a supplemented base for $\psi_m(EU_{2n}(R_m,$ $\Lambda_m))$, there is a $U'_1 \in A'$ such that $^{d'_1}U_1'\subseteq U'$. Clearly
\begin{align*}
&^{U'}(^{d'_1}g')\\
\supseteq~~&^{U'}(^{d'_1} (^{U_1'}g'))\\
\overset{(4.9.1)}{\supseteq}& ^{U'}(^{d'_1} V') \text{ (for some } V' \in B'
 \text{)}\\
\overset{(4.9.2)}{\supseteq}&^{U'}g''  \text{ (for some nontrivial image } g'' \text{ of an elementary matrix in } EU_{2n}(R,\Lambda)\text{)} \\
\overset{(4.9.1)}{\supseteq}&V'_1 \text{ (for some } V'_1 \in B' \text{).}
\end{align*}
\underline{case 2} Assume that $k > 1$. Set $h' := d'_{k-1}\dots d'_1$. Thus $d' = d'_k\dots d'_1 = d'_kh'$.  We can assume by induction 
on $k$ that given $U'_1 \in A'$, $^{U'_1}(^{h'}g')\supseteq V'$ for some $V' \in B'$. Now we proceed similarly to case 1, replacing $g'$ by $^{h'}g'$ and $d'_1$ by $d'_k$. Here are the details. Choose $U'_1\in A'$ such that $^{d'_k}\phi(U'_1)\subseteq U'$. Clearly
\begin{align*}
&^{U'}(^{d'_kh'}g')\\
\supseteq~~&^{U'}(^{d'_k} (^{U_1'}(^{h'}g')))\\
\overset{I.A.}{\supseteq}~&^{U'}(^{d'_k} V')\\
\overset{(4.9.2)}{\supseteq}&^{U'}g''  \text{ (for some nontrivial image } g'' \text{ of an elementary matrix in } EU_{2n}(R,\Lambda)\text{)} \\
\overset{(4.9.1)}{\supseteq}&V'_1 \text{ (for some } V'_1 \in B' \text{).}
\end{align*}
\hfill$\qed$
}
\Theorem{Let $n$ and $(R,\Lambda)$ be as in the third paragraph of this section. Let $H$ be a subgroup of $U_{2n}(R,\Lambda)$. Then 
\begin{align*} &H \text{ is normalized by } EU_{2n}(R,\Lambda)\Leftrightarrow \\ 
 &\exists!\text{ form ideal }(I,\Gamma)\text{ such that }EU_{2n}((R,\Lambda),(I,\Gamma))\subseteq H\subseteq CU_{2n}((R,\Lambda),(I,\Gamma)).
 \end{align*}
}
\Proof{$~$\\
$\Rightarrow$:\\
Assume that $H$ is normalized by $EU_{2n}(R,\Lambda)$. We have to show existence and uniqueness of a form ideal $(I,\Gamma)$ such that 
\[EU_{2n}((R,\Lambda),(I,\Gamma))\subseteq H\subseteq CU_{2n}((R,\Lambda),(I,\Gamma)).\]
\underline{existence}\\
Set $I:=\{x\in R|T_{12}(x)\in H\}$ and $\Gamma:=\{y\in \Lambda|T_{-1,1}(y)\in H\}$. Then $(I,\Gamma)$ is a form ideal, $EU_{2n}((R,\Lambda),(I,\Gamma))\subseteq H$ and $(I,\Gamma)$ is maximal with this property (i.e. if $EU_{2n}((R,\Lambda),(I',\Gamma'))\subseteq H$, then $I'\subseteq I$ and $\Gamma'\subseteq \Gamma$). It remains to show that $H\subseteq CU_{2n}((R,\Lambda),(I,\Gamma))$. The proof is by contradiction. Suppose $H\not\subseteq CU_{2n}((R,\Lambda),(I,\Gamma))$. Then the image $\hat H$ of $H$ in $U_{2n}(R,\Lambda)/U_{2n}((R,\Lambda),(I,\Gamma))$ contains a noncentral element $\hat h$ by the definition of $CU_{2n}((R,\Lambda),(I,\Gamma))$. By Lemma 4.1 there is a maximal ideal $m$ of $C$ such that $I\cap C\subseteq m$ and $h':=\phi_m(\hat h)$ is noncentral in $U_{2n}(R_m,\Lambda_m)/U_{2n}((R_m,\Lambda_m),(I_m,\Gamma_m))$. Choose an $s_0\in S_m$ with the properties (1) and (2) in Lemma 4.2, let $(A,B)$ be the supplemented base for $EU_{2n}(R,\Lambda)$ defined in Lemma 4.4 and set $(A',B'):=\phi_m(\psi(A,B))$. Choose an $U'\in A'$. By Lemma 4.8 there is a $k\in\mathbb{N}$ and elements $g'_0,\dots,g'_k\in\phi_m(\psi(EU_{2n}(R,\Lambda)))$, $\epsilon'_0,\dots,\epsilon'_k\in \psi_m(EU_{2n}(R_m,$  $\Lambda_m))$ and $l_1,\dots,l_k\in\{-1,1\}$ such that $g'_k$ is the nontrivial image of an elementary matrix in $EU_{2n}(R,\Lambda)$, 
\[ ^{\epsilon'_{k}}([^{\epsilon'_{k-1}}(\dots^{\epsilon'_2}([^{\epsilon'_1}([^{\epsilon'_0}h',g'_0]^{l_1}),g'_1]^{l_2})\dots),g'_{k-1}]^{l_{k}})=g'_k\tag{4.11.1}\]
and
\[^{d'_i}g'_i\in U'~\forall i\in \{0,\dots,k\}\]
where $d'_i=(\epsilon'_{i}\cdot\hdots\cdot\epsilon'_0)^{-1}\hspace{0.1cm}\forall i\in\{0,\dots,k\}$. By conjugating (4.11.1) by $d'_k=(\epsilon'_{k}\cdot\hdots\cdot\epsilon'_0)^{-1}$ we get
\[[\dots[[h',^{d'_0}\hspace{-0.1cm}g'_0]^{l_1},^{d'_1}\hspace{-0.1cm}g'_1]^{l_2}\dots,^{d'_{k-1}}\hspace{-0.1cm}g'_{k-1}]^{l_{k}}=^{d'_k}\hspace{-0.1cm}g'_k.\tag{4.11.2}\]
By Corollary 4.10 there is a $V'\in B'$ such that 
\[V'\subseteq {}^{U'}(^{d'_k}g'_k).\tag{4.11.3}\]
Let $\hat U\in \hat A:=\psi(A)$ and $\hat V\in \hat B:=\psi(B)$ such that $\phi_m(\hat U)=U'$ and $\phi_m(\hat V)=V'$. Clearly we may assume that $\hat V\subseteq \hat U$ since $\hat V\cap \hat U$ contains a member of $\hat B$. Since $^{d'_i}g'_i\in U'=\phi_m(\hat U)~\forall i\in\{0,\dots,k\}$, there are $x_0,\dots,x_{k}\in \hat U$ such that $\phi_m(x_i)={}^{d'_i}g'_i~\forall i\in\{0,\dots,k\}$. Set 
\[x:=[\dots[[\hat h,x_0]^{l_1},x_1]^{l_2}\dots,x_{k-1}]^{l_k}.\]
Clearly the l.h.s. of $(4.11.2)$ equals $\phi_m(x)$. Further $x\in \hat H$ since $\hat H$ is normalized by $\psi(EU_{2n}(R,\Lambda))$ and $x\in M:=\psi(U_{2n}((R,\Lambda),
(s_0R,s_0\Lambda)))$  since $x_{k-1}\in \hat U\subseteq M$ and $M$ is normal. It follows that $^{\hat U}x\subseteq \hat H\cap M$. By (4.11.2) and (4.11.3) we have 
\[\phi_m(\hat V)\subseteq{}^{\phi_m(\hat U)}\phi_m(x)=\phi_m({}^{\hat U}x)\subseteq\phi_m(\hat H\cap M).\tag{4.11.4}\]
By Lemma 4.2, $\phi_m$ is injective on $M$. Hence it follows from (4.11.4) that $\hat V\subseteq \hat H$ (note that $\hat V\subseteq \hat U\subseteq M$). Let $V=EU_{2n}(Rxs_0R,\Gamma(xs_0))\in B$, where $x\in R, xs_0\not\in I$ or $x\in\Lambda, xs_0\not\in \Gamma$, such that $\hat V=\psi(V)$. It follows that $V\subseteq H\cdot U_{2n}((R,\Lambda),(I,\Gamma))$. This implies $EU_{2n}((R,\Lambda),(Rxs_0R,\Gamma(xs_0)))\subseteq H\cdot U_{2n}((R,\Lambda),(I,\Gamma))$, since both $H$ and $U_{2n}((R,\Lambda),(I,\Gamma))$ are normalized by $EU_{2n}(R,\Lambda)$. Hence 
\begin{align*}
&EU_{2n}((R,\Lambda),(Rxs_0R,\Gamma(xs_0)))\\
=&[EU_{2n}(R,\Lambda),EU_{2n}((R,\Lambda),(Rxs_0R,\Gamma(xs_0)))]\\
\subseteq&[EU_{2n}(R,\Lambda),H\cdot U_{2n}((R,\Lambda),(I,\Gamma))]\\
\subseteq&[EU_{2n}(R,\Lambda),H]({}^H[EU_{2n}(R,\Lambda),U_{2n}((R,\Lambda),(I,\Gamma))])\\
=&[EU_{2n}(R,\Lambda),H]({}^HEU_{2n}((R,\Lambda),(I,\Gamma)))\\
\subseteq& H
\end{align*}
by Lemma 3.16. But this contradicts the maximality of $(I,\Gamma)$ since clearly $Rxs_0R\not\subseteq I$ or $\Gamma(xs_0)\not\subseteq \Gamma$. Thus $H\subseteq CU_{2n}((R,\Lambda),(I,\Gamma))$.\\
\\
\underline{uniqueness}\\
Assume that \[EU_{2n}((R,\Lambda),(I,\Gamma))\subseteq H\subseteq CU_{2n}((R,\Lambda),(I,\Gamma))\]
and
\[EU_{2n}((R,\Lambda),(I',\Gamma'))\subseteq H\subseteq CU_{2n}((R,\Lambda),(I',\Gamma')).\]
It follows that 
\begin{align*}
&EU_{2n}((R,\Lambda),(I,\Gamma))\\
=&[EU_{2n}(R,\Lambda),EU_{2n}((R,\Lambda),(I,\Gamma))]\\
\subseteq& [EU_{2n}(R,\Lambda),CU_{2n}((R,\Lambda),(I',\Gamma'))]\\
=&EU_{2n}((R,\Lambda),(I',\Gamma')).
\end{align*}
It is easy to deduce that $I\subseteq I'$ and $\Gamma\subseteq \Gamma'$. By symmetry it follows that $I=I'$ and $\Gamma=\Gamma'$.\\
\\
$\Leftarrow$:\\
Suppose that \[EU_{2n}((R,\Lambda),(I,\Gamma))\subseteq H\subseteq CU_{2n}((R,\Lambda),(I,\Gamma)).\]
Then 
\begin{align*}
[EU_{2n}(R,\Lambda),H]
\subseteq [EU_{2n}(R,\Lambda),CU_{2n}((R,\Lambda),(I,\Gamma))]
=EU_{2n}((R,\Lambda),(I,\Gamma))\subseteq H
\end{align*}
and hence $H$ is normalized by $EU_{2n}(R,\Lambda)$.
\hbox{}\hfill$\qed$}
\Definition{Let $R$ be a ring, $r\mapsto \bar r$ an involution on $R$ and $\lambda\in center(R)$. Then we call $R$ {\it quasi-finite} if it is a direct limit of subrings $R_i$ ($i\in \Phi$ where $\Phi$ is some index set) which are almost commutative, involution invariant and contain $\lambda$ (recall that a ring is called almost commutative if it is module finite over its center). 
} 
\Lemma[A. Bak]{Let $T$ be an almost commutative ring, $t\mapsto \bar t$ an involution on $T$ and $\lambda\in center(T)$. Then $T$ is a direct limit of involution invariant subrings $T_{j}~(j\in \Psi)$ containing $\lambda$ such that for any $j\in \Psi$, $T_{j}$ is a Noetherian $C_j$-module where $C_j$ is the subring of $T_j$ consisting of all finite sums of elements of the form $c\bar c$ and $-c\bar c$ where $c$ ranges over some subring $C'_j\subseteq center(T_j)$.
}
\Proof{Denote the center of $T$ by $C$. Since $T$ is almost commutative, there is an $q\in\mathbb{N}$ and elements $x_1,\dots,x_q\in T$ such that $T = Cx_1+\dots+Cx_q$. For each $k,l\in\{1,\dots,q\}$ there are $a^{(kl)}_1,\dots, a^{(kl)}_q\in C$ such that $x_kx_l=\sum\limits_{p=1}^{q}a^{(kl)}_px_p$. Further for each $k\in\{1,\dots,q\}$ there are $b^{(k)}_1,\dots, b^{(k)}_q\in C$ such that $\bar x_k=\sum\limits_{p=1}^{q}b^{(k)}_px_p$. Finally there are $c_1,\dots, c_q\in C$ such that $\lambda =\sum\limits_{p=1}^{q}c_px_p$. Set \[K:=\mathbb{Z}[a^{(kl)}_p,\overline{a^{(kl)}_p},b^{(k)}_p,\overline{b^{(k)}_p},c_p,\overline{c_p}|k,l,p\in\{1,\dots,q\}].\]
One checks easily that $C$ is a $K$-algebra and the direct limit of all involution invariant $K$-subalgebras $A_{j}~(j\in \Psi)$ of $C$ which are finitely generated over $K$.
For any $j\in \Psi$ set $T_{j}:=A_{j}+A_{j}x_1+\dots+A_{j}x_q$. One checks easily that each $T_{j}$ is an involution invariant subring of $T$ containing $\lambda$. Further $\varinjlim\limits_{j} T_{j}=T$. Fix a $j\in \Psi$ and let $C_{j}$ denote the subring of $A_{j}$ consisting of all finite sums of elements of the form $a\bar a$ and $-a\bar a$ where $a\in A_{j}$. We will show that $T_j$ is a Noetherian $C_j$-module. Clearly $A_{j}$ is a finitely generated $\mathbb{Z}$-algebra and hence also a finitely generated $C_{j}$-algebra. Since for any $a\in A_{j}$
\[a + \bar a = (a + 1)(\bar a + 1) - a\bar a-1,\]
$C_{j}$ contains all sums $a+\bar a$ where $a\in A_{j}$. Since any $a\in A_{j}$ is root of the monic polynomial $X^2-(a+\bar a)X+a\bar a$, $A_{j}$ is an integral extension of $C_{j}$. Since $A_{j}$ is an integral extension of $C_{j}$ and a finitely generated $C_{j}$-algebra, $A_{j}$ is a finitely generated module over $C_{j}$ by \cite{lang}, chapter VII, Proposition 1.2. Since $T_{j}$ is finitely generated over $A_{j}$, it is a finitely generated $C_{j}$-module. Since $K$ is a Noetherian ring, $A_{j}$ is a Noetherian ring (by Hilbert's Basis Theorem) and hence $C_{j}$ is a Noetherian ring (by the Eakin-Nagata Theorem). Thus $T_{j}$ is a Noetherian $C_{j}$-module.\hfill $\qed$
}
\Theorem{Let $n\geq 3$ and $(R,\Lambda)$ a form ring where $R$ is quasi-finite. Let $H$ be a subgroup of $U_{2n}(R,\Lambda)$. Then 
\begin{align*} &H \text{ is normalized by } EU_{2n}(R,\Lambda)\Leftrightarrow \\ 
 &\exists!\text{ form ideal }(I,\Gamma)\text{ such that }EU_{2n}((R,\Lambda),(I,\Gamma))\subseteq H\subseteq CU_{2n}((R,\Lambda),(I,\Gamma)).
 \end{align*}
}
\Proof{Let $(I,\Gamma)$ denote the level of $H$, i.e. the largest form ideal such that $EU_{2n}((R,\Lambda),(I,\Gamma))\subseteq H$. We will show that $H \subseteq CU_{2n}((R,\Lambda),(I,\Gamma))$, i.e. $[\sigma,\tau]\in U_{2n}((R,\Lambda),(I,\Gamma))$ for any  $\sigma\in H$ and $\tau\in U_{2n}(R,\Lambda)$. Let $\sigma\in H$ and $\tau\in U_{2n}(R,\Lambda)$. Since $R$ is quasi-finite, it is the direct limit of almost commutative, involution invariant subrings $R_i~(i\in \Phi)$ containing $\lambda$. By the the previous lemma each $R_i$ is the direct limit of involution invariant subrings $R_{ij}~(j\in \Psi_i)$ containing $\lambda$ such that for any $j\in \Psi_i$, $R_{ij}$ is a Noetherian $C_{ij}$-module where $C_{ij}$ is the subring of $R_{ij}$ consisting of all finite sums of elements of the form $c\bar c$ and $-c\bar c$ where $c$ ranges over some subring $C'_{ij}\subseteq center(R_{ij})$. If $i\in \Phi$ and $j\in \Psi_i$, set $\Lambda_{ij}:=\Lambda\cap R_{ij}$. One checks easily that $(R_{ij},\Lambda_{ij})$ is a form ring. Clearly there is an $i\in \Phi$ and a $j\in \Psi_i$ such that $\sigma,\tau\in U_{2n}(R_{ij},\Lambda_{ij})$. Set $H_{ij}:=U_{2n}(R_{ij},\Lambda_{ij})\cap H$. Then $\sigma\in H_{ij}$ and $H_{ij}$ is normalized by $EU_{2n}(R_{ij},\Lambda_{ij})$. Let $(I_{ij},\Gamma_{ij})$ denote the level of $H_{ij}$. Then obviously $I_{ij}\subseteq I$ and $\Gamma_{ij}\subseteq\Gamma$. By Theorem 4.11,
\[H_{ij} \subseteq CU_{2n}((R_{ij},\Lambda_{ij}),(I_{ij},\Gamma_{ij})).\]
Hence $[\sigma,\tau]\in U_{2n}((R_{ij},\Lambda_{ij}),(I_{ij},\Gamma_{ij}))
\subseteq U_{2n}((R,\Lambda),(I,\Gamma))$. Thus we have shown that 
\[EU_{2n}((R,\Lambda),(I,\Gamma))\subseteq H\subseteq CU_{2n}((R,\Lambda),(I,\Gamma))\]
where $(I,\Gamma)$ is the level of $H$. The uniqueness of $(I,\Gamma)$ and the implication $\Leftarrow$ follow from the standard commutator formulas (see the proof of Theorem 4.11).\hfill $\qed$}\\

\end{document}